\documentclass[final,1p,times]{elsarticle}

\usepackage{amssymb}
\usepackage{bm}
 \usepackage{amsthm}
\usepackage{amscd}
\usepackage{amsmath}
\usepackage{amsfonts}
\usepackage{amssymb}
\usepackage{graphicx}
\usepackage{color}
\newtheorem{theorem}{Theorem}

\newtheorem{lemma}[theorem]{Lemma}

\usepackage{mathrsfs}
\usepackage{titletoc}

\newcommand{\ra}{\rightarrow}
\newcommand{\p}{\partial}
\newcommand{\f}{\frac}

\newcommand{\be}{\begin{equation}}
\renewcommand{\ra}{\rightarrow}
\newcommand{\ee}{\end{equation}}
\newcommand{\bea}{\begin{eqnarray}}
\newcommand{\eea}{\end{eqnarray}}
\newcommand{\bna}{\begin{eqnarray*}}
\newcommand{\ena}{\end{eqnarray*}}

\renewcommand{\le}{\left}
\newcommand{\ri}{\right}

\journal{***}

\begin{document}

\begin{frontmatter}
\title{Blow-up analysis involving isothermal coordinates on the boundary of compact Riemann surface}
\author{Yunyan Yang}
 \ead{yunyanyang@ruc.edu.cn}
 \address{Department of Mathematics,
Renmin University of China, Beijing 100872, China}

\author{Jie Zhou}
\ead{zhoujie2014@mails.ucas.ac.cn}
\address{Department of Mathematics, Tsinghua University, Beijing 100084, China}

\begin{abstract}
  Using the method of blow-up analysis, we obtain two sharp Trudinger-Moser inequalities on a compact Riemann surface with smooth boundary,
  as well as the existence of the corresponding extremals.
  This generalizes early results of Chang-Yang \cite{Chang-Yang} and the first named author \cite{2006}, and complements Fontana's inequality
  of two dimensions \cite{Fontana}.
  The blow-up analysis in the current paper is far more elaborate than that of \cite{2006},
  and particularly clarifies several ambiguous points there. In precise,
  we prove the existence of isothermal coordinate systems near the boundary, the existence and uniform estimates of the Green function with the Neumann boundary condition. Also our analysis can be applied to the Kazdan-Warner problem and the Chern-Simons Higgs problem on compact Riemman surfaces with smooth boundaries.
\end{abstract}


\end{frontmatter}

\section{Introduction}

Let $\Omega$ be a smooth bounded domain in $\mathbb{R}^2$, $W_0^{1,2}(\Omega)$ be the completion of all smooth functions with compact
 support under the norm
 $$\|u\|_{W_0^{1,2}(\Omega)}=\le(\int_\Omega |\nabla u|^2dx\ri)^{1/2}.$$
 It was proved by Yudovich \cite{Yudovich}, Pohozaev \cite{Pohozaev}, Peetre \cite{Peetre}, Trudinger \cite{Trudinger} and Moser \cite{Moser}
 that
\be\label{Moser}\sup_{u\in W_0^{1,2}(\Omega),\,\|u\|_{W_0^{1,2}(\Omega)}\leq 1}\int_\Omega\exp(\gamma u^2)dx<+\infty,\quad\forall\gamma\leq 4\pi;\ee
 moreover, if $\gamma>4\pi$, then the above supremum is infinity. In literature,
 such kind of inequalities are known as Trudinger-Moser inequalities.
 Concerning all smooth functions with mean value zero instead of boundary value zero, Chang and Yang \cite{Chang-Yang} obtained
 by using their isoperimetric inequality that
\be\label{C-Y}\sup_{u\in W^{1,2}(\Omega),\,\int_\Omega|\nabla u|^2dx\leq 1,\,\int_\Omega udx=0}\int_\Omega\exp(\gamma u^2)dx<+\infty,\quad \forall\gamma\leq 2\pi.\ee
Analogous to (\ref{Moser}), the supremum in (\ref{C-Y}) is infinity for any $\gamma>2\pi$. This inequality was applied  by Chang and Yang to the Nirenberg problem with the Neumann boundary
condition.

Now we consider $(\Sigma,g)$, a closed Riemann surface, i.e. a compact Riemann surface without boundary, and let $W^{1,2}(\Sigma, g)$ be the usual Sobolev space. Representing a function
by the Riesz potential of its gradient and using a manifold version of Adams' potential estimate \cite{Adams}, L. Fontana was able to show the following:
\be\label{Fontana}\sup_{u\in W^{1,2}(\Sigma,g),\,\int_\Sigma|\nabla_gu|^2dv_g\leq 1,\,\int_\Sigma udv_g=0}\int_\Sigma\exp(\gamma u^2)dv_g<+\infty,\quad
\forall\gamma\leq 4\pi,\ee
where $\nabla_g$ and $dv_g$ stand for the gradient operator and the Riemann volume element with respect to the metric $g$. Similar to the Euclidean case, the above supremum
is infinity for any $\gamma>4\pi$. Of course, L. Fontana obtained far more than (\ref{Fontana}) in his elegant paper \cite{Fontana}. Later, via a method of blow-up analysis, Li \cite{LiJPDE} proved that the supremum in (\ref{Fontana}) can be attained for all $\gamma\leq 4\pi$.

In view of (\ref{C-Y}), one would naturally expect (\ref{Fontana}) for compact Riemann surfaces with smooth boundaries. Indeed, in the case that $(\Sigma,g)$ is a compact Riemann surface with smooth boundary $\p\Sigma$, following the approach of Li \cite{LiJPDE}, the first named author \cite{2006} extended  (\ref{C-Y}) as below:
\bea\label{06-1}\sup_{u\in W^{1,2}(\Sigma,g),\,\int_\Sigma|\nabla_gu|^2dv_g\leq 1,\,\int_\Sigma udv_g=0}\int_\Sigma\exp(\gamma u^2)dv_g<+\infty,&\gamma\leq 2\pi,\\
\label{06-2}\sup_{u\in W^{1,2}(\Sigma,g),\,\int_\Sigma(|\nabla_gu|^2+u^2)dv_g\leq 1}\int_\Sigma\exp(\gamma u^2)dv_g<+\infty,&\gamma\leq 2\pi.\eea
Furthermore, both supremums can be attained for all $\gamma\leq 2\pi$, but they are infinite when $\gamma>2\pi$.

Let us revisit the outline of the proof of (\ref{06-1}) \cite{2006} (the proof of (\ref{06-2}) is almost the same as that of (\ref{06-1})). For any $k\in\mathbb{N}$, by a direct method of variation, there exists a $u_k\in W^{1,2}(\Sigma,g)\cap C^1(\overline{\Sigma})$
with $\int_\Sigma|\nabla_gu_k|^2dv_g=1$ and $\int_\Sigma u_kdv_g=0$ such that $u_k$ is a maximizer for the supremum in (\ref{06-1}) with  $\gamma=\gamma_k=2\pi-1/k$, and $u_k$ satisfies the
Euler-Lagrange equation
\be\label{E-L-0}\le\{
\begin{array}{lll}
\Delta_gu_k=\f{1}{\lambda_k}u_k\exp(\gamma_ku_k^2)-\f{\mu_k}{\lambda_k}&{\rm in}&\Sigma\\[1.2ex]
\p u_k/\p\bm\nu=0&{\rm on}& \p\Sigma\\[1.2ex]
\lambda_k=\int_\Sigma u_k^2\exp(\gamma_ku_k^2)dv_g\\[1.2ex]
\mu_k=\int_\Sigma u_k\exp(\gamma_ku_k^2)dv_g,
\end{array}
\ri.\ee
where $\Delta_g$ is the Laplace-Beltrami operator and $\bm\nu$ denotes the unit outward vector fields on $\p\Sigma$. With no loss of generality, we can assume
$c_k=u_k(x_k)=\max_\Sigma |u_k|\ra+\infty$ and $x_k\ra x_0\in\p\Sigma$ as $k\ra\infty$. To proceed,  we {\it firstly}
 choose an isothermal coordinate system $(U,\psi;\{y_1,y_2\})$ near $x_0$ satisfying $\psi(x_0)=(0,0)$,
$\psi(\p U\cap\p\Sigma)\subset\p\mathbb{R}^{2+}$, and the metric $g$ can be written as $g=\exp(2f(y))(dy_1^2+dy_2^2)$, where
$f\in C^1(\overline{U})$ and $f(0,0)=0$. Choosing appropriate scale $r_k>0$ and applying elliptic estimates, we have up to a subsequence,
 $$\label{bubble-0}c_k\le(u_k\circ\psi^{-1}(\psi({x}_k)+r_k\cdot)-c_k\ri)\ra -\f{1}{2\pi}\log\le(1+\f{\pi}{2}|\cdot|^2\ri)\quad{\rm in}\quad
 C^1_{\rm loc}(\mathbb{R}^{2+}\cup\p\mathbb{R}^{2+}).$$
 {\it Secondly} we  prove that $c_ku_k$ converges to some Green function $G_{x_0}$ weakly in $W^{1,q}(\Sigma,g)$ for any $1<q<2$, strongly
 in $L^s(\Sigma,g)$ with $s<2q/(2-q)$, and in $C^1_{\rm loc}(\overline{\Sigma}\setminus\{x_0\})$, where $G_{x_0}$ is a distributional solution
 of the equation
 $$\le\{\begin{array}{lll}
 \Delta_gG_{x_0}=\delta_{x_0}-\f{1}{{\rm Area}(\Sigma)}&{\rm in}&\Sigma\\[1.2ex]
 \p G_{x_0}/\p\bm\nu=0&{\rm on}&\p\Sigma\\[1.2ex]
 \int_\Sigma G_{x_0}dv_g=0.
 \end{array}\ri.$$
 Based on elliptic estimates in the isothermal coordinate system near $x_0$, $G_{x_0}$ can be locally decomposed as
 \be\label{Green-x0}G_{x_0}(x)=-\f{1}{\pi}\log{\rm dist}_g(x_0,x)+A_{x_0}+O({\rm dist}_g(x_0,x)).\ee
 {\it Thirdly}, using the capacity estimate introduced by Li \cite{LiJPDE}, we derive
 \be\label{upbd-x0}\sup_{u\in W^{1,2}(\Sigma,g),\,\int_\Sigma|\nabla_gu|^2dv_g\leq 1,\,\int_\Sigma udv_g=0}
 \int_\Sigma\exp(2\pi u^2)dv_g\leq {\rm Area}(\Sigma)+\f{\pi}{2}\exp(1+2\pi A_{x_0}).\ee
 {\it Finally} we construct a function sequence $\phi_k\in W^{1,2}(\Sigma,g)$ with $\int_\Sigma\phi_kdv_g=0$
 and $\int_\Sigma|\nabla_g\phi_k|^2dv_g=1$ satisfying
 \be\label{test-x0}\int_\Sigma\exp(2\pi \phi_k^2)dv_g> {\rm Area}(\Sigma)+\f{\pi}{2}\exp(1+2\pi A_{x_0}),\ee
 provided that $k$ is chosen sufficiently large. The contradiction between (\ref{upbd-x0}) and (\ref{test-x0}) implies
 that $c_k$ must be bounded. Then applying elliptic estimates to (\ref{E-L-0}), one has up to a subsequence, $u_k\ra u_0$ in $C^1(\overline{\Sigma})$
 as $k\ra\infty$ and thus the supremum in (\ref{06-1}) can be attained by $u_0$.

 Checking the proof in \cite{2006}, we found at least three key points that should have been seriously treated there.
 The first one is the claimed existence of isothermal coordinate system on the boundary $\p\Sigma$, which is very important in the
 subsequent blow-up analysis;
 The second one is the way of finding a constant $C$ depending only on $(\Sigma,g)$ and $q<2$ such that
 \be\label{bdd0}\int_\Sigma|\nabla_g(c_ku_k)|^qdv_g\leq C,\ee
 which leads to the convergence of $c_ku_k$ to $G_{x_0}$; The third one is the decomposition of $G_{x_0}$ with the form (\ref{Green-x0}).

 Our goals are twofold. One is to clarify the above three concerns. Specifically, we employ Riemann mapping theorems to construct isothermal coordinate
 systems near the boundary $\p\Sigma$; To prove (\ref{bdd0}), we first construct a Green function with the Neumann boundary condition, and then use the Green representation formula; The decomposition of $G_{x_0}$ will be based on
 elliptic estimates in an isothermal coordinate system. The other one is to improve Theorems 1.1 and 1.2 in \cite{2006}.
 To describe this improvement, we define a space of functions by
\be\label{H-space}\mathcal{H}=\le\{u\in W^{1,2}(\Sigma,g): \int_\Sigma udv_g=0\ri\},\ee
the first eigenvalue of the Laplace-Beltrami operator with respect to the Neumann boundary condition by
\be\label{eigen}\lambda_{\rm N}(\Sigma)=\inf_{u\in \mathcal{H},\,u\not\equiv 0}\f{\int_\Sigma|\nabla u|^2dv_g}{\int_\Sigma u^2dv_g},\ee
 and a Sobolev norm on $\mathcal{H}$ in the case $\alpha<\lambda_{\rm N}(\Sigma)$ by
\be\label{1-alpha}\|u\|_{1,\alpha}=\le(\int_\Sigma|\nabla u|^2dv_g-\alpha\int_\Sigma u^2dv_g\ri)^{1/2}.\ee

 Our main result reads as follows:
 \begin{theorem}\label{thm1}
 Let $(\Sigma,g)$ be a compact Riemann surface with smooth boundary $\p\Sigma$. Then for any $\alpha<\lambda_{\rm N}(\Sigma)$, there holds
 \be\label{ineq-1}\sup_{u\in\mathcal{H},\,\|u\|_{1,\alpha}\leq 1}\int_\Sigma\exp(\gamma u^2)dv_g<+\infty,\quad\forall\gamma\leq 2\pi,\ee
 where $\mathcal{H}$, $\lambda_{\rm N}(\Sigma)$ and $\|\cdot\|_{1,\alpha}$ are defined as in (\ref{H-space}), (\ref{eigen}) and (\ref{1-alpha})
 respectively. Moreover, the above supremum is infinity for any $\gamma>2\pi$.  Furthermore,  for any fixed $\alpha<\lambda_{\rm N}(\Sigma)$
 and $\gamma\leq 2\pi$, the supremum in (\ref{ineq-1}) can be attained
 by some $u^\ast\in\mathcal{H}\cap C^1(\overline{\Sigma})$ with
 $\|u^\ast\|_{1,\alpha}=1$.
 \end{theorem}

 Similarly we have the following:

 \begin{theorem}\label{thm2}
 Let $(\Sigma,g)$ be a compact Riemann surface with smooth boundary $\p\Sigma$. Then for any real number $\tau>0$, there holds
 \be\label{ineq-2}\sup_{u\in W^{1,2}(\Sigma,g),\,\|u\|_{1,\tau}\leq 1}\int_\Sigma\exp(\gamma u^2)dv_g<+\infty,
 \quad\forall\gamma\leq 2\pi,\ee
 where $\|u\|_{1,\tau}=(\int_\Sigma(|\nabla_gu|^2+\tau u^2)dv_g)^{1/2}$.
 Moreover, the above supremum is infinity for any $\gamma>2\pi$. Furthermore,  for all real numbers $\tau>0$ and $\gamma\leq 2\pi$, the supremum in (\ref{ineq-2}) can be attained by some $u_0\in\mathcal{H}\cap C^1(\overline{\Sigma})$ with
 $\|u_0\|_{1,\tau}=1$.
 \end{theorem}

 Theorems \ref{thm1} and \ref{thm2} are complements of \cite{Chang-Yang,Fontana,LiJPDE,Yang-Tran,Yang-JDE-15}.
 We remark that the inequality (\ref{ineq-1}) involving the norm $\|\cdot\|_{1,\alpha}$ was motivated by \cite{Tintarev}, while
 the inequality (\ref{ineq-2}) involving the norm $\|\cdot\|_{1,\tau}$ was motivated by Adimurthi-Yang \cite{Adi-Yang} and do \'O-Yang
 \cite{do-Yang}.
 Although the method of blow-up analysis is now standard,
 the technique is far more delicate than the existing related works \cite{DJLW,Adi-Stru,LiJPDE,2006}. Our technique can certainly be used in
  the study of Trudinger-Moser inequalities on boundaries
 \cite{Liu,Li-Liu,Yang-PJM,Yang-MZ,Lu-Yang,Nguyen}, as well as in the Chern-Simons Higgs problem with Neumann boundary condition \cite{DJLW1,DJLW2,Lan-Li,WangM1,WangM-2},
 and other related problems \cite{Bonheure,Deng-Musso,Deng,ZhouChunqin,ZhuX}.

 As far as the inequality itself is concerned, (\ref{ineq-2}) is apparently weaker than (\ref{ineq-1}), but unexpectedly
 they are equivalent. Motivated by \cite{Nguyen-2}, we have
 the following:

 \begin{theorem}\label{thm3}
 Let $(\Sigma,g)$ be a compact Riemann surface with smooth boundary $\p\Sigma$, and $\lambda_{\rm N}(\Sigma)$ be defined as in
 (\ref{eigen}). Given any $0\leq\alpha<\lambda_{\rm N}(\Sigma)$ and any $\tau>0$. Assume that (\ref{ineq-2})
 holds for all $\gamma<2\pi$. Then the inequality
 \be\label{in-01}\sup_{u\in\mathcal{H},\,\|u\|_{1,\alpha}\leq 1}\int_\Sigma\exp(2\pi u^2)dv_g<+\infty\ee
 is equivalent to
 \be\label{in-02}\sup_{u\in W^{1,2}(\Sigma,g),\,\|u\|_{1,\tau}\leq 1}\int_\Sigma\exp(2\pi u^2)dv_g<+\infty,
 \ee
 where $\mathcal{H}$, $\|\cdot\|_{1,\alpha}$ and $\|\cdot\|_{1,\tau}$ are the same as in Theorems \ref{thm1} and \ref{thm2} respectively.
 \end{theorem}

 Throughout this paper, sequence and subsequence are not distinguished, and various constants are often denoted by the same $C$.
 The remaining part of this paper is organized as follows: In Section \ref{Sec-2}, we prove the existence of isothermal coordinate system around any point
 on the boundary $\p\Sigma$; In Section \ref{Sec-3}, we construct a Green function with the Neumann boundary condition and give its
 uniform estimates; Theorems \ref{thm1}-\ref{thm3} will be proved in Sections \ref{Sec-4}-\ref{Sec-6} respectively.

\section{Isothermal coordinate systems near the boundary}\label{Sec-2}

In this section, we prove existence of isothermal coordinate systems near the boundary. This is based on the classical existence result near  inner points
of Riemann surface and  Riemann mapping theorems involving the boundary. From now on, we always denote
$$\mathbb{B}_r^+=\le\{y=(y_1,y_2)\in\mathbb{R}^2: y_1^2+y_2^2<r,\,\,y_2>0\ri\},\,\,
\mathbb{R}^{2+}=\le\{y=(y_1,y_2)\in\mathbb{R}^2: y_2\geq 0\ri\}$$
and the closure of a set $E$ by $\overline{E}$.

\begin{lemma}\label{prop}
Let $(\Sigma,g)$ be a compact Riemann surface with smooth boundary $\p\Sigma$. For any fixed point $x\in\p\Sigma$,
there exist a number $\delta>0$ and
an isothermal coordinate system $(\overline{U}_x,\psi_x;\{y_1,y_2\})$ near $x$ such that $\psi_x(x)=(0,0)$, $\overline{U}_x\subset\overline{\Sigma}$ is
a neighborhood of $x$,
$\psi_x(U_x)=\mathbb{B}_\delta^+$ and $\psi_x(\overline{U}_x\cap\p\Sigma)=\overline{\mathbb{B}_\delta^+}\cap \p{\mathbb{R}^{2+}}$.
In this coordinate system, there exists a function $f\in C^1(\overline{\mathbb{B}_\delta^+}, \mathbb{R})$ such that for all $y=(y_1,y_2)\in \overline{\mathbb{B}_\delta^+}$, the metric $g$ can be written as
$$\label{metric}g=\exp{(2f(y))}(dy_1^2+dy_2^2).$$ Suppose that $\bm\nu$ is an unit outward vector field defined on $\psi_x^{-1}(\overline{\mathbb{B}_\delta^+}\cap \p{\mathbb{R}^2}^+)\subset\p\Sigma$. For any $p\in\psi_x^{-1}(\overline{\mathbb{B}_\delta^+}\cap \p{\mathbb{R}^2}^+)$, if we write $y=\psi_x(p)$, then
$$\label{outward}(\psi_x)_\ast(\bm\nu(p))=\exp(-f(y)){\p}/{\p y_2}.$$
\end{lemma}

{\it Proof}.
We divide the construction into several steps.

{\it Step 1. There exists a neighborhood $\overline{U}_1$ of $x$, a domain $\Omega_1\subset\mathbb{R}^2$ verifying that $\p\Omega_1$
  is smooth except for two corners, and a homeomorphism $\psi_1:\overline{U}_1\ra\overline{\Omega}_1$ such that $\psi_1(x)=(0,0)$ and
$\psi_1(\overline{U}_1\cap \p\Sigma)=\Gamma_1\subset\p\Omega_1$. In the coordinate system $(\overline{U}_1,\psi_1;\{x_1,x_2\})$, the metric
$g$ can be written as $g=\exp(2f_1(x_1,x_2))(dx_1^2+dx_2^2)$ for all $(x_1,x_2)\in\overline{\Omega}_1$, where $f_1$ is a smooth function
with $f(0,0)=0$. Denote $\bm\nu_1=(\psi_1)_\ast(\bm\nu)$.
Then $\bm\nu_1=\exp(-f_1(x_1,x_2))\nu_0$, where $\bm\nu_0$ is the unit outward vector field on $\p\Omega_1$.}

Indeed, since $(\Sigma,g)$ is a compact
Riemann surface with smooth boundary $\p\Sigma$,
we understand that there exists another compact Riemann surface $(\Sigma^\ast,{g}^\ast)$ with smooth boundary $\p\Sigma^\ast$
such that $\overline{\Sigma}\subset {\Sigma}^\ast$, ${\rm dist}_g(\overline{\Sigma},\p {\Sigma}^\ast)>0$ and ${g}^\ast=g$ on
$\overline{\Sigma}$. Note that $x$ is an inner point of $\Sigma^\ast$.
By \cite{Bers}, there exist $U\subset{\Sigma}^\ast$, a neighborhood of $x$,  and
a diffeomorphism $\psi_1: U\ra \mathbb{B}_r\subset \mathbb{R}^2$ with $\psi_1(x)=(0,0)$ such that the metric $g$ reads as
\be\label{metr-local}g=\exp{(2f_1(x_1,x_2))}(dx_1^2+dx_2^2),\ee
where $f_1$ is a smooth function with $f_1(0,0)=0$.
Denote $U_1=U\cap\Sigma$ and $\Omega_1=\psi_1(U_1)$. To finish this step, it suffices to estimate $\bm\nu_1$.
Write $\bm\nu_1=a_1\p/\p x_1+a_2\p/\p x_2$. Then
\be\label{equality}1=|\bm\nu|^2=\exp(2f_1(x_1,x_2))(a_1^2+a_2^2),\ee
which immediately leads to the representation of $\bm\nu_1$.

{\it Step 2}. Replace $\Omega_1$ by a smooth domain ${\Omega}_2\subset\Omega_1$ verifying that $\psi_1(x)=(0,0)$ is an inner point of
a smooth curve ${\Gamma}_2\subset \Gamma_1\cap \p{\Omega}_2$.

{\it Step 3}. {\it $\overline{\Omega}_2$ is conformal to a unit disc $\overline{\mathbb{D}}\subset\mathbb{R}^2$.}
In fact, according to the Riemann mapping theorem \cite{Bers}, there exists a conformal map $\psi_2:{\Omega}_2\ra\mathbb{D}$
denoted by $w=\psi_2(z)$ with $\psi_2(0,0)=(0,-1)$, where $z=x_1+ix_2$.
By (\cite{Pom}, Theorem 3.5), $\psi_2$ extends to a map in $C^1(\overline{\Omega}_2,\overline{\mathbb{D}})$; moreover,
$\psi_2^\prime(z)\not=0$ for all $z\in \overline{\Omega}_2$. Here and in the sequel we slightly abuse some notations.
In particular we
identify $z\in\mathbb{C}$ with $(x_1,x_2)\in\mathbb{R}^2$, and so on.

{\it Step 4}. {\it $\mathbb{D}$ is conformal to a half plane.}
 Let $q\not\in \psi_2(\Gamma_2)$ be fixed. Then via a M\"obius transformation $\zeta=h(w)$,
the set $\mathbb{D}\setminus\{q\}$ can be mapped
into the upper half plane $\mathbb{R}^{2+}$ with $h(\psi_2(\psi_1(x)))=(0,0)$.
Define a function by
\be\label{phi}\varphi(z)=f_1(z)-\log|h^\prime(w)\psi_2^\prime(z)|\ee
and a dilation $\tau: \mathbb{R}^{2+}\ra\mathbb{R}^{2+}$ by
$y=\tau(\zeta)=\exp(\varphi(0,0))\zeta$.
Thus
\be\label{dy}dy=\exp(\varphi(0,0))h^\prime(w)\psi_2^\prime(z)dz.\ee
Set $\psi_x=\tau\circ h\circ \psi_2\circ\psi_1$.
Choose $\delta>0$ sufficiently small so that $\psi_x^{-1}(\p\mathbb{B}_\delta^+
 \cap \p\mathbb{R}^{2+})\subset \psi_1^{-1}(\Gamma_2)$.

{\it Step 5}. {\it $(\psi_x^{-1}(\overline{\mathbb{B}_\delta^+}), \psi_x;\{y_1,y_2\})$ is
an isothermal coordinate system near $x\in\p\Sigma$ as we required.} Indeed, since $\psi_1$, $\psi_2$, $h$
and $\tau$ are all conformal maps, we conclude that $\psi_x$ is also a conformal map. This together with (\ref{metr-local}),
(\ref{phi}) and (\ref{dy})
leads to the representation of
the metric $g$ as
\bna
g&=&\f{\exp{(2f_1(z))}}{|h^\prime(w)\psi_2^\prime(z)|^2}\exp{(-2\varphi(0,0))}|dy|^2\\
&=&\exp{(2\varphi(z)-2\varphi(0,0))}|dy|^2\\
&=&\exp{(2f(y))}(dy_1^2+dy_2^2)
\ena
for all $y\in \overline{\mathbb{B}_\delta^+}$, where
$f(y)=\varphi(z)-\varphi(0,0)$ and $z=\psi_2^{-1}(h^{-1}(\exp(-\varphi(0,0))y))$. By the above definitions of
$h$ and $\psi_2$, we have that $y=(0,0)$
if and only if $z=(0,0)$. Thus $f(0,0)=0$.
 Moreover, we can assume $(\psi_x)_\ast(\bm\nu)(p)=b_2(y)\p/\p y_2$
for any $p\in \psi_x^{-1}(\mathbb{B}_\delta^+)\cap\p\Sigma$ with $y=\psi_x(p)$. Similar to (\ref{equality}),
we calculate $b_2(y)=\exp(-f(y))$. Clearly $f\in C^0(\overline{\mathbb{B}_\delta^+})$. Further application of (\cite{Pom}, Theorem 3.6)
implies that $f$ is smooth on $\overline{\mathbb{B}_\delta^+}$.
This ends the proof of the lemma. $\hfill\Box$

\section{The Green function with the Neumann boundary condition}\label{Sec-3}

In this section, we concern the Green function on $(\Sigma,g)$
with the Neumann boundary condition, whose construction is based on the method of (Aubin \cite{Aubin}, Chapter 4).
For its uniform estimate, we use elliptic estimates as Aubin did in (\cite{Aubin}, Chapter 4), and as Druet, Robert, Wei did in \cite{DRW}.
To begin with, we need the following:

\begin{lemma}\label{W22}
  Let $(\Sigma,g)$ be a compact Riemann surface with smooth boundary $\p\Sigma$. If $f\in L^2(\Sigma,g)$ satisfies $\int_\Sigma fdv_g=0$,
  then there exists a unique weak solution of
  \be\label{E-L-N}\le\{\begin{array}{lll}\Delta_g u=f&{\rm in}&\Sigma\\[1.2ex]
  {\p u}/{\p\bm\nu}=0&{\rm on}&\p\Sigma\\[1.2ex]
  \int_\Sigma u dv_g=0,\end{array}\ri.
  \ee
  or equivalently there exists a $u\in\mathcal{H}$ defined by (\ref{H-space}) satisfies
  \be\label{weak-sol}\int_\Sigma\nabla_g u\nabla_g\varphi dv_g=\int_\Sigma f\varphi dv_g,\quad \forall \varphi\in C^1(\overline{\Sigma}).\ee
  Moreover there exists some constant $C$ depending only on $(\Sigma,g)$ such that
  \be\label{W2-est}
  \|u\|_{W^{2,2}(\Sigma,g)}\leq C\|f\|_{L^2(\Sigma,g)}.
  \ee
   If further $f\in C^\alpha(\overline{\Sigma})$ for some $0<\alpha<1$, then $u\in C^{2,\alpha}(\overline{\Sigma})$.
  \end{lemma}

  {\it Proof}. The uniqueness is obvious. To see this, we let $u_1$ and $u_2$ be two weak solutions of (\ref{E-L-N}) and $u^\ast=u_1-u_2$.
  Since $C^1(\overline{\Sigma})$ is dense in $W^{1,2}(\Sigma,g)$, it follows from (\ref{weak-sol}) that
  $$\int_\Sigma\nabla_gu^\ast\nabla_g vdv_g=0,\quad\forall v\in W^{1,2}(\Sigma,g).$$
  Choosing $v=u^\ast$ in the above equality, we conclude $u^\ast\equiv 0$ since $u^\ast\in\mathcal{H}$.

  The Existence of weak solution of (\ref{E-L-N}) is based on a direct method of variation. Let us consider the functional
  $$J(u)=\f{1}{2}\int_\Sigma |\nabla_gu|^2dv_g-\int_\Sigma fudv_g.$$
  For any $u\in\mathcal{H}$, we have by the H\"older inequality and the Poincare inequality
  $$\label{Holder}\int_\Sigma fudv_g\leq C\le(\int_\Sigma f^2dv_g\ri)^{1/2}\le(\int_\Sigma|\nabla_gu|^2dv_g\ri)^{1/2},$$
  which implies that $J$ has a lower bound on $\mathcal{H}$. Now we
  take a sequence of functions $u_j\in\mathcal{H}$ satisfying $J(u_j)\ra \inf_{u\in\mathcal{H}}J(u)$. One can easily see that
  $u_j$ is bounded in $\mathcal{H}$. Thus one can assume up to a subsequence, $u_j$ converges to some $u_0\in \mathcal{H}$
  weakly in $W^{1,2}(\Sigma,g)$, strongly in $L^q(\Sigma,g)$ for any $q>1$ and almost everywhere in $\Sigma$. Clearly
  $u_0\in\mathcal{H}$ and
  $$J(u_0)\leq \lim_{j\ra\infty}J(u_j)=\inf_{u\in\mathcal{H}}J(u).$$
  Hence $u_0$ is a minimizer of $J$ on $\mathcal{H}$ and satisfies the Euler-Lagrange equation (\ref{weak-sol}).

  We now prove (\ref{W2-est}). Since $(\Sigma,g)$ is compact, we have by the standard $W^{2,2}$-estimate (see for example \cite{Aubin}, Theorem 3.54) that
  \be\label{H2}\|u_0\|_{W^{2,2}(\Sigma,g)}\leq C(\|u_0\|_{L^2(\Sigma,g)}+\|f\|_{L^2(\Sigma,g)}^2)\ee
  for some constant $C$ depending only on $(\Sigma,g)$. Noting that $u_0$ is a unique solution of (\ref{E-L-N}), we have by the definition of
   distributional solution and the H\"older inequality that
  $$\int_\Sigma |\nabla_g u_0|^2dv_g=\int_\Sigma fu_0dv_g\leq \|u_0\|_{L^2(\Sigma,g)}\|f\|_{L^2(\Sigma,g)}.$$
  This together with the Poincare inequality leads to
  \be\label{Poinc}\|u_0\|_{L^2(\Sigma,g)}\leq C\|f\|_{L^2(\Sigma,g)}.\ee
  Inserting  (\ref{Poinc}) into (\ref{H2}), we conclude (\ref{W2-est}), as desired.

  Finally, if $f\in C^\alpha(\overline{\Sigma})$, then we have $u\in C^{2,\alpha}(\overline{\Sigma})$ by using Lemma \ref{prop} and
  the classical Schauder estimate (\cite{GT}, Theorem 6.6).
  $\hfill\Box$\\

An analog of (\cite{Aubin}, Theorems 4.13 and 4.17) reads as follows.
\begin{lemma}\label{G-N}
There exists a unique Green function ${\textsf{G}}(x,\cdot)\in L^1(\Sigma,g)$ satisfying
  \be\label{Green-equation}\le\{\begin{array}{lll}
  \Delta_{g,y}{\textsf{G}}(x,y)=\delta_x(y)-\f{1}{{\rm Area}(\Sigma)}&{\rm in}&\Sigma\\[1.5ex]
  \f{\p}{\p \nu_y}\textsf{G}(x,y)=0&{\rm on}&\p\Sigma\\[1.5ex]
  \int_\Sigma \textsf{G}(x,y)dv_{g,y}=0
  \end{array}\ri.\ee
  in the distributional sense,
  or equivalently for any $\varphi\in C^2(\overline{\Sigma})$ with $\p\varphi/\p\bm\nu=0$
  on $\p\Sigma$, there holds
  \be\label{G-weak}\int_{\Sigma}{\textsf{G}(x,y)}\Delta_g\varphi(y)dv_{g,y}=\varphi(x)-\overline{\varphi},\ee
  where $\overline{\varphi}=\f{1}{{\rm Area}(\Sigma)}\int_\Sigma\varphi dv_g$.
  Moreover, for all $x,y\in \overline{\Sigma}$ with $x\not= y$, $\textsf{G}(x,y)=\textsf{G}(y,x)$,
  and there exists some constant $C$ depending only on $(\Sigma,g)$ such that
  \be
  |\textsf{G}(x,y)|\leq C(1+|\log{\rm dist}_g(x,y)|),\quad|\nabla_{g,y}\textsf{G}(x,y)|\leq {C}({\rm dist}_g(x,y))^{-1},
  \label{grad}\ee
  where ${\rm dist}_g(x,y)$ denotes the geodesic distance between $x$ and $y$.
  \end{lemma}
{\it Proof}.
   {\it Part I. Uniqueness of the Green function}. If $\textsf{G}_1(x,y)$ and $\textsf{G}_2(x,y)$ are two Green functions satisfying (\ref{G-weak}),
   then we set $h(y)=\textsf{G}_1(x,y)-\textsf{G}_2(x,y)$. By Lemma \ref{W22}, for any $f\in C^\alpha(\overline{\Sigma})$, $0<\alpha<1$, there exists a unique
   $\varphi\in C^{2,\alpha}(\overline{\Sigma})$ such that $\Delta_g\varphi=f-\overline{f}$, $\p\varphi/\p\bm\nu=0$ on $\p\Sigma$, and
   $\int_\Sigma \varphi dv_g=0$. Hence (\ref{G-weak}) implies that
   $$\int_\Sigma h fdv_g=\int_\Sigma h\Delta_g\varphi dv_g=0.$$
   This together with the facts $h\in L^1(\Sigma,g)$ and $C^\alpha(\overline{\Sigma})$ is dense in $L^1(\Sigma,g)$ leads to $h\equiv 0$.

   {\it Part II. Existence of the Green function}.

  {\it Case 1. $x$ is an inner point of $\Sigma$.} We follow the line of (\cite{Aubin}, Theorem 4.13). Let $i(x)$ be the injectivity radius of $x$,  and $\phi(r)$ be a decreasing function, which is
  equal to $1$ in a neighborhood of zero, and to zero for $r>i(x)/8$. Define
  $$H(x,y)=-\f{1}{2\pi}\phi({\rm dist}_g(x,y))\log {\rm dist}_g(x,y),$$
  $\Gamma(x,y)=\Gamma_1(x,y)=\Delta_{g,y}H(x,y)$, $\Gamma_{i+1}(x,y)=\int_\Sigma\Gamma_i(x,z)\Gamma(z,y)dv_{g,z}$ for $i=1, 2$,
  and set
  $$\label{GF}\textsf{G}(x,y)=H(x,y)+\sum_{i=1}^2\int_\Sigma\Gamma_i(x,z)H(z,y)dv_{g,z}+F(x,y),$$
  where $F(x,y)$ satisfies
  $$\label{Fxy}\le\{\begin{array}{lll}\Delta_{g,y}F(x,y)=-\Gamma_{3}(x,y)-\f{1}{{\rm Area}(\Sigma)}&{\rm in}&\Sigma\\[1.2ex]
  \f{\p}{\p\bm\nu_y}F(x,y)=0&{\rm on}& \p\Sigma\end{array}
  \ri.$$
  and $$\int_\Sigma F(x,y)dv_{g,y}=-\int_\Sigma H(x,y)dv_{g,y}-\sum_{j=1}^2\int_\Sigma\le(\int_\Sigma\Gamma_i(x,z)H(z,y)
  dv_{g,z}\ri)dv_{g,y}.$$
  Such an $F(x,y)$ exists in view of Lemma \ref{W22}. It can be easily checked that $\textsf{G}(x,\cdot)$ satisfies (\ref{G-weak}).

  {\it Case 2. $x\in \p\Sigma$.} By Lemma \ref{prop}, we choose an isothermal coordinate system $(\overline{U_x},\psi_x;\{z_1,z_2\})$ near $x$ such that
  $\psi_x: U_x\ra\mathbb{B}_{r_0}^+$ for some $r_0>0$. In this coordinate system, the metric $g$ can be written as
  $$g=\exp(2f(z))(dz_1^2+dz_2^2);$$
  moreover $\Delta_g=-\exp(-2f(z))\Delta_{\mathbb{R}^2}$. Let $\phi:[0,\infty)\ra[0,\infty)$ be a
  decreasing function such that $\phi\equiv 1$ on $[0,r_0/2]$ and $\phi\equiv 0$ on $[r_0,\infty)$. Set
  $$H(x,y)=\le\{\begin{array}{lll}
  -\f{1}{\pi}\phi(|z|)\log|z|,&y=\psi_x^{-1}(z)\in\psi_x^{-1}\le(\overline{\mathbb{B}_{r_0}^+}\ri)\\[1.2ex]
  0,&y\in\Sigma\setminus\psi_x^{-1}\le(\overline{\mathbb{B}_{r_0}^+}\ri).
  \end{array}\ri.$$
  One can check that
  $$\le\{\begin{array}{lll}
  \Delta_{g,y}H(x,y)=\delta_x(y)-\eta(x,y)&{\rm in}&\Sigma\\[1.2ex]
  \f{\p}{\p\bm\nu_y}H(x,y)=0&{\rm on}& \p\Sigma
  \end{array}\ri.$$
  in the distributional sense (\ref{G-weak}), where
  $$\eta(x,y)=\le\{\begin{array}{lll}
  \exp(-2f(z))\Delta_{\mathbb{R}^2}\le(-\f{1}{\pi}\phi(|z|)\log|z|\ri), &y=\psi_{x}^{-1}(z)\in\psi_x^{-1}\le(\overline{\mathbb{B}_{r_0}^+}\ri)\\[1.2ex]
  0,&y\in\Sigma\setminus\psi_x^{-1}\le(\overline{\mathbb{B}_{r_0}^+}\ri).
  \end{array}\ri.$$
  According to Lemma \ref{W22}, one can find a unique $F(x,y)$ satisfying
  $$\le\{\begin{array}{lll}
  \Delta_{g,y}F(x,y)=\eta(x,y)-\f{1}{{\rm Area}(\Sigma)}&{\rm in}&\Sigma\\[1.2ex]
  \f{\p}{\p\bm\nu_y}F(x,y)=0&{\rm on}& \p\Sigma\\[1.2ex]
  \int_\Sigma F(x,y)dv_{g,y}=-\int_\Sigma H(x,y)dv_{g,y}.
  \end{array}\ri.$$
  We set
  $\textsf{G}(x,y)=H(x,y)+F(x,y)$ for all $y\in\overline{\Sigma}$. Then $\textsf{G}(x,\cdot)$ is a distributional solution of
  (\ref{Green-equation}).

  {\it Part III. Uniform estimate}.

  We first prove that there exists some constant $C$ depending only on $(\Sigma,g)$ such that for all $x\in\overline{\Sigma}$,
  there holds
  \be\label{Lp}\|\textsf{G}(x,\cdot)\|_{L^2(\Sigma,g)}\leq C.\ee
  To see this, for any $w\in C^2(\overline{\Sigma})$, we conclude from Lemma \ref{W22} that the equation
  $$\label{equ-0}\le\{\begin{array}{lll}\Delta_g u=w-\f{1}{{\rm Area}(\Sigma)}\int_\Sigma wdv_g&{\rm in}&\Sigma\\[1.2ex]
  {\p u}/{\p\bm\nu}=0&{\rm on}&\p\Sigma\\[1.2ex]
  \int_\Sigma udv_g=0\end{array}\ri.
  $$
    has a unique solution  $u\in C^2(\overline{\Sigma})$.
    Combining (\ref{Green-equation}) and (\ref{W2-est}), we obtain
    $$
    \int_\Sigma \textsf{G}(x,y)w(y)dv_{g,y}=\int_\Sigma \textsf{G}(x,y)\Delta_gu(y)dv_{g,y}
    \leq\|u\|_{C^0(\overline{\Sigma})}\leq C\|u\|_{W^{2,2}(\Sigma,g)}\leq C\|w\|_{L^2(\Sigma,g)},
    $$
    where $C$ is a constant depending only on $(\Sigma,g)$.
    This together with the density of $C^2(\overline{\Sigma})$ in $L^2(\Sigma,g)$ implies (\ref{Lp}).

   Given any fixed $x_0\in\p\Sigma$. Take an isothermal coordinate system $(\overline{U_{x_0}},\psi_{x_0};\{z_1,z_2\})$ near $x_0$ such that
   $\psi(x_0)=(0,0)$, $\psi_{x_0}(\overline{U_{x_0}})=\overline{\mathbb{B}_{r_0}^+}$ and $\psi_{x_0}^{-1}(\overline{\mathbb{B}_{r_0}^+}\cap\p\mathbb{R}^{2+})=\overline{U_{x_0}}\cap\p\Sigma$
   for some $r_0>0$, and the metric $g=\exp(2f(z))(dz_1^2+dz_2^2)$ with $f\in C^2(\overline{\mathbb{B}_{r_0}^+})$ and $f(0)=0$. For any point $x\in \psi_{x_0}^{-1}(\mathbb{B}_{r_0/8}^+)$, we define
   $${G}_x^\ast(z)=\le\{\begin{array}{lll}
   \textsf{G}(x,\psi_{x_0}^{-1}(z_1,z_2)),&z=(z_1,z_2)\in\overline{\mathbb{B}_{r_0}^+}\\[1.2ex]
   \textsf{G}(x,\psi_{x_0}^{-1}(z_1,-z_2)),&z=(z_1,z_2)\in\mathbb{B}_{r_0}\setminus\overline{\mathbb{B}_{r_0}^+}.
   \end{array}\ri.$$
   Then ${G}_{x}^\ast$ is a distributional solution of
   \be\label{distr}-\Delta_{\mathbb{R}^2}{G}_x^\ast(z)=\delta_{z_0}(z)+\delta_{z_0^\prime}(z)-\f{\exp(2f^\ast(z))}{{\rm Area}(\Sigma)}\quad
   {\rm in}\quad\mathbb{B}_{r_0},\ee
   where $z_0=\psi_{x_0}(x)=(z_{0,1},z_{0,2})$, $z_0^\prime=(z_{0,1},-z_{0,2})$ and
   $$f^\ast(z)=\le\{\begin{array}{lll}
   f(z_1,z_2),&z_2\geq 0\\[1.2ex]
   f(z_1,-z_2),&z_2<0.
   \end{array}\ri.$$
   Denote
   \be\label{FG}{F}_x^\ast(z)={G}_x^\ast(z)+\f{1}{2\pi}\log|z-z_0|+\f{1}{2\pi}\log|z-z_0^\prime|.\ee
   It follows from (\ref{distr}) that ${F}_x^\ast$ satisfies
   \be\label{F}\Delta_{\mathbb{R}^2}{F}_x^\ast(z)=\f{\exp(2f^\ast(z))}{{\rm Area}(\Sigma)}\quad{\rm in}\quad\mathbb{B}_{r_0}\ee
   in the distributional sense. By (\ref{Lp}), we have
   $\|{F}_x^\ast(\cdot)\|_{L^2(\mathbb{B}_{r_0})}\leq C$ for some constant $C$ depending only on $(\Sigma,g)$ and $r_0$.
   Then applying $W^{2,2}$-estimate to (\ref{F}), we can see that ${F}_x^\ast$ is bounded in $W^{2,2}(\mathbb{B}_{2r_0/3})$
   uniformly with respect to $x\in\psi_{x_0}^{-1}(\overline{\mathbb{B}_{r_0/2}^+})$. Further elliptic estimate leads to
   $$\|{F}_x^\ast(\cdot)\|_{C^1(\overline{\mathbb{B}_{r_0/4}})}\leq C$$
   for some constant $C$ depending only on $(\Sigma,g)$ and $r_0$. This together with (\ref{FG}) gives
   $$|{G}_x^\ast(z)|\leq C(1+\log|z-z_0|+\log|z-z_{0}^\prime|)\leq C(1+\log|z-z_0|)$$
   and
   $$|\nabla_{\mathbb{R}^2}{G}_x^\ast(z)|\leq C(|z-z_0|^{-1}+|z-z_0^\prime|^{-1})\leq C|z-z_0|^{-1}$$
   for all $z\in\overline{\mathbb{B}_{r_0/4}^+}\setminus\{z_0\}$, since $|z-z_0|<|z-z_0^\prime|$.
   Therefore there exists some constant $C$ depending only on $(\Sigma,g)$ and $r_0$ such that
   \be\label{C1-1}
   |\textsf{G}(x,y)|\leq C(1+|\log {\rm dist}_g(x,y)|),\,\,
   |\nabla_g\textsf{G}(x,y)|\leq C({\rm dist}_g(x,y))^{-1}
   \ee
   for all $x\in\psi_{x_0}^{-1}(\overline{\mathbb{B}_{r_0/8}^+})$ and $y\in\psi_{x_0}^{-1}(\overline{\mathbb{B}_{r_0/4}^+)}$ with $y\not=x$.
   Now for any fixed $x\in\psi_{x_0}^{-1}(\overline{\mathbb{B}_{r_0/8}})$, in view of (\ref{Lp}), we have by applying
   elliptic estimate to (\ref{Green-equation}) that
   $$\|\textsf{G}(x,\cdot)\|_{C^1(\overline{\Sigma\setminus\psi_{x_0}^{-1}(\mathbb{B}_{r_0/4})})}\leq C$$
   for some constant $C$ depending only on $(\Sigma,g)$ and $r_0$. This implies (\ref{C1-1}) already holds for
   all $x\in\psi_{x_0}^{-1}(\overline{\mathbb{B}_{r_0/8}^+})$ and all $y\in\Sigma$ with $y\not=x$.
   Since $\p\Sigma$ is compact, one can find a real number $r_1>0$ and a constant $C$
   depending only on $(\Sigma,g)$ and $r_1$ such that
   \be\label{conv-bdr}(\ref{C1-1})\,\,{\rm holds\,\,for \,\,all}\,\, x\in \Sigma_{r_1}=\le\{x\in\overline{\Sigma}: {\rm dist}_g(x,\p\Sigma)\leq r_1\ri\} \,\,{\rm and}\,\,y\in\overline{\Sigma}\,\,
   {\rm with}\,\, y\not=x.\ee

   If $x_0$ is an inner point of $\Sigma$, we take an isothermal coordinate system $(U_{x_0},\psi_{x_0};\{z_1,z_2\})$ near $x_0$ such
   that $\psi_{x_0}(x_0)=(0,0)$, $U_{x_0}\subset \Sigma\setminus\Sigma_{r_1/2}$,  $\psi_{x_0}(U_{x_0})=\mathbb{B}_{r_0}$,
   and the metric
   $g=\exp(2f(z))(dz_1^2+dz_2^2)$ with $f(0,0)=0$. For any $x\in \psi_{x_0}^{-1}(\mathbb{B}_{r_0/8})$,
    we define $G_x(z)=\textsf{G}(x,\psi_{x_0}^{-1}(z))$ for $z\in\mathbb{B}_{r_0}$.
   Denote $z_0=\psi_{x_0}(x)=(z_{0,1},z_{0,2})$. Then $G_x(z)$ is a distributional solution of
   $$-\Delta_{\mathbb{R}^2}G_x(z)=\delta_{z_0}(z)-\f{\exp(2f(z))}{{\rm Area}(\Sigma)}\quad{\rm in}\quad \mathbb{B}_{r_0}.$$
   As a consequence
   \be\label{decom}-\Delta_{\mathbb{R}^2}\le({G}_x(z)+\f{1}{2\pi}\log|z|\ri)=-\f{\exp(2{f}(z))}{{\rm Area}(\Sigma)}\quad{\rm in}\quad
   \mathbb{B}_{r_0}\ee
   in the distributional sense.
   In view of (\ref{Lp}), applying elliptic estimate to (\ref{decom}), we conclude ${G}_x(z)+\f{1}{\pi}\log|z|$ is bounded
   in $C^1(\overline{\mathbb{B}_{r_{0}/2}})$ uniformly in $x\in\psi_{x_0}^{-1}(\mathbb{B}_{r_0/8})$, and thus
   \be\label{local-2}\le\|{\textsf{G}}(x,\cdot)+\f{1}{2\pi}\log{\rm dist}_g(x,\cdot)\ri\|_{C^1(\overline{\psi_{x_0}^{-1}(\mathbb{B}_{r_0/2}}))}\leq C\ee
   for all $x\in\psi_{x_0}^{-1}(\mathbb{B}_{r_0/8})$, where $C$ is a constant depending only on $(\Sigma,g)$, $\psi_{x_0}$ and $r_{0}$.
   In addition, we have by applying elliptic estimate to (\ref{Green-equation}) that
   $\|\textsf{G}(x,\cdot)\|_{C^1(\overline{\Sigma\setminus\psi_{x_0}^{-1}(\mathbb{B}_{r_0/2})})}\leq C$ for all $x\in\psi_{x_0}^{-1}(\mathbb{B}_{r_0/8})$.
   This together with (\ref{local-2}) implies that (\ref{C1-1}) holds for some constant $C$ depending only on $(\Sigma,g)$, $\psi_{x_0}$ and $r_{0}$,
   and for all $x\in\psi_{x_0}^{-1}(\mathbb{B}_{r_0/8})$. As a result, in view of the compactness of $\overline{{\Sigma}\setminus\Sigma_{r_1}}$, we conclude that there exists some constant $C$, depending only on $(\Sigma,g)$ and $r_1$, such that
   \be\label{conv-bdr-2}(\ref{C1-1})\,\,{\rm holds\,\,for \,\,all}\,\, x\in \overline{\Sigma\setminus\Sigma_{r_1}} \,\,{\rm and}\,\,y\in\overline{\Sigma}\,\,
   {\rm with}\,\, y\not=x.\ee

   Combining (\ref{conv-bdr}) and (\ref{conv-bdr-2}), we conclude (\ref{grad}), as desired.

   {\it Part IV. Symmetry}.
   {We shall prove that $\textsf{G}(x,y)=\textsf{G}(y,x)$ for all $(x,y)\in\overline{\Sigma}\times\overline{\Sigma}$ with $x\not=y$.}

   For any $f\in C^\alpha(\overline{\Sigma})$, we set
   \be\label{F-x}F(x)=\int_\Sigma \textsf{G}(y,x)(f(y)-\overline{f})dv_{g,y},\ee
   where $\overline{f}=\f{1}{{\rm Area}(\Sigma)}\int_\Sigma fdv_g$. In view of
   (\ref{grad}), we have $\textsf{G}\in L^1(\Sigma\times\Sigma)$. Hence we obtain by the Fubini theorem
   \be\label{F-0}\overline{F}=\f{1}{{\rm Area}(\Sigma)}\int_\Sigma \le(\int_\Sigma\textsf{G}(y,x)dv_{g,x}\ri)(f(y)-\overline{f})dv_{g,y}=0.\ee
   By Lemma \ref{W22}, there exists a unique $\varphi\in C^{2,\alpha}(\overline{\Sigma})$ satisfying
   \be\label{varphi}\le\{\begin{array}{lll}\Delta_g \varphi=f-\overline{f}&{\rm in}&\Sigma\\[1.2ex]
  {\p \varphi}/{\p\bm\nu}=0&{\rm on}&\p\Sigma\\[1.2ex]
  \overline{\varphi}=0.\end{array}\ri.\ee
   We now {\it claim} that
   \be\label{F=var}F(x)=\varphi(x)\,\,\,{\rm for\,\,all}\,\,\, x\in \overline{\Sigma}.\ee
   By (\ref{grad}) and the Lebesgue dominated convergence theorem, one can easily see that
   $F$ is continuous on $\overline{\Sigma}$. For any $h\in C^\alpha(\overline{\Sigma})$, there exists a unique
   $\psi\in C^{2,\alpha}(\overline{\Sigma})$ such that
   \be\label{psi}\le\{\begin{array}{lll}\Delta_g \psi=h-\overline{h}&{\rm in}&\Sigma\\[1.2ex]
  {\p \psi}/{\p\bm\nu}=0&{\rm on}&\p\Sigma\\[1.2ex]
  \overline{\psi}=0.\end{array}\ri.\ee
  By (\ref{F-x}), (\ref{F-0}), (\ref{varphi}), (\ref{psi}) and the Fubini theorem, we calculate
  \bna
  \int_\Sigma Fhdv_{g}&=&\int_\Sigma F(h-\overline{h})dv_{g}\\
  &=&\int_\Sigma\le(\int_\Sigma\textsf{G}(y,x)(f(y)-\overline{f})dv_{g,y}\ri)\Delta_g\psi(x)dv_{g,x}\\
  &=&\int_\Sigma\le(\int_\Sigma\textsf{G}(y,x)\Delta_g\psi(x)dv_{g,x}\ri)(f(y)-\overline{f})dv_{g,y}\\
  &=&\int_\Sigma\psi(y)\Delta_{g,y}\varphi(y)dv_{g,y}\\
  &=&\int_\Sigma(h-\overline{h})\varphi(y)dv_{g,y}\\
  &=&\int_\Sigma h\varphi dv_g.
  \ena
  Noting that $h\in C^\alpha(\overline{\Sigma})$ is arbitrary, $F\in C^0(\overline{\Sigma})$ and $\varphi\in C^{2,\alpha}(\overline{\Sigma})$,
  we conclude (\ref{F=var}).

  It follows from (\ref{G-weak}) and (\ref{F=var}) that
  \bna
  \int_\Sigma  \textsf{G}(y,x)(f(y)-\overline{f})dv_{g,y}=\int_\Sigma  \textsf{G}(x,y)\Delta_g\varphi(y)dv_{g,y}
  =\int_\Sigma  \textsf{G}(x,y)(f(y)-\overline{f})dv_{g,y}.
  \ena
  As a consequence
  \be\label{G0}\int_\Sigma(\textsf{G}(y,x)-\textsf{G}(x,y))(f(y)-\overline{f})dv_{g,y}=0.\ee
  Denote $\mu(x)=\f{1}{{\rm Area}(\Sigma)}\int_\Sigma\textsf{G}(y,x)dv_{g,y}$. Clearly $\mu\in C^0(\overline{\Sigma})$
  because of (\ref{grad}). Since $f\in C^\alpha(\overline{\Sigma})$ is arbitrary, we conclude from (\ref{G0}),
  $$\textsf{G}(y,x)-\textsf{G}(x,y)=\mu(x)\quad{\rm for\,\,a. e.}\,\, y\in\Sigma.$$
  Integrating both sides of the above equation with respect to $x$, we have by the Fubini theorem
  $$\mu(y)=\f{1}{{\rm Area}(\Sigma)}\int_\Sigma\textsf{G}(x,y)dv_{g,x}=-\f{1}{{\rm Area}(\Sigma)}\int_\Sigma\mu(x)dv_{g,x}=0,$$
  which implies that $\mu\equiv 0$ on $\overline{\Sigma}$, and whence
  $\textsf{G}(y,x)=\textsf{G}(x,y)$  for a.e. $(x,y)\in\Sigma\times\Sigma$. Since $\textsf{G}(x,\cdot)\in C^1(\overline{\Sigma}\setminus\{x\})$
  due to (\ref{grad}), we have $\textsf{G}(\cdot,x)\in C^1(\overline{\Sigma}\setminus\{x\})$. Therefore $\textsf{G}(x,y)$ is continuous
  for all $(x,y)\in\overline{\Sigma}\times\overline{\Sigma}$ with $x\not=y$, and this gives the symmetry of $\textsf{G}(\cdot,\cdot)$.
    $\hfill\Box$

\section{Proof of Theorem \ref{thm1}}\label{Sec-4}
In this section, we shall prove Theorem \ref{thm1} by using the method of blow-up analysis. Pioneer woks related to this topic are
due to Ding-Jost-Li-wang \cite{DJLW}, Adimurthi-Struwe \cite{Adi-Stru}, and Li \cite{LiJPDE}. Here, in our situation, blow-up
happens on the boundary $\p\Sigma$. This brings new difficulties compared with the previous situation \cite{LiJPDE,Yang-JDE-15}.
In particular, we use the Green representation formula of $c_ku_k$ to obtain the boundedness of $\|\nabla_g(c_ku_k)\|_{L^q(\Sigma,g)}$
for any $1<q<2$, which is the key step in the study of the convergence of $c_ku_k$ (see Lemma \ref{Green} below). It should be mentioned
that our blow-up analysis and decomposition of certain Green function depend on the
existence of isothermal coordinate system near the boundary $\p\Sigma$.

Since the proof is very long, we sketch it as follows:  In the subsection \ref{4.1},  let $\gamma^\ast$ be the best constant for the inequality (\ref{ineq-1}),
  which will be explicitly defined by (\ref{gamma}) below. Then $\gamma^\ast$ must be $2\pi$.  In the subsections \ref{4.2} and \ref{4.3}, there exists a smooth maximizer for any subcritical
  Trudinger-Moser functional. If blow-up happens (the maximizers are not uniformly bounded), by a process of blow-up analysis on a sequence of maximizers, we obtain an accurate estimate on the supremum in
   (\ref{ineq-1}). In the subsection \ref{4.4}, we construct a sequence of admissible functions to show that the supremum in (\ref{ineq-1}) is strictly greater
   than that we obtained in the subsection \ref{4.3}. This implies that no blow-up happens in the subsection \ref{4.3}. Thus  elliptic estimate leads to the attainability of
   the supremum in (\ref{ineq-1}) for $\gamma=2\pi$.
\subsection{The best constant}\label{4.1}

Let $\mathcal{H}$, $\lambda_{\rm N}(\Sigma)$ and $\|\cdot\|_{1,\alpha}$ be defined as in (\ref{H-space})-(\ref{1-alpha})
respectively. We first have

\begin{lemma}\label{lemma1}
For any $\alpha<\lambda_{\rm N}(\Sigma)$, there exists some constant $\gamma_0>0$ such that
$$\sup_{u\in \mathscr{H},\,\|u\|_{1,\alpha}\leq 1}\int_\Sigma \exp(\gamma_0u^2)dv_g<+\infty.$$
\end{lemma}
{\it Proof}.
Since $\alpha<\lambda_{\rm N}(\Sigma)$, we have for any $u$ with $\|u\|_{1,\alpha}\leq 1$,
\be\label{eng}
1\geq\int_\Sigma|\nabla_gu|^2dv_g-\alpha\int_\Sigma u^2dv_g
\geq\le(1-\f{\alpha}{\lambda_{N}(\Sigma)}\ri)\int_\Sigma|\nabla_gu|^2dv_g.\ee
If $x$ is an inner point of $\Sigma$, we choose an isothermal coordinate system $(U_x,\psi_x;\{y_1,y_2\})$ around $x$,
where $U_x\subset\Sigma$ is a neighborhood of $x$ and $\psi_x: U_x\ra\Omega_x\subset\mathbb{R}^2$ is a diffeomorphism. In this
coordinate system,
the metric $g=\exp(2f(y))(dy_1^2+dy_2^2)$, where $f\in C^1(\overline{\Omega_x})$. As a consequence, we have by (\ref{eng}) that
\be\label{euc}\int_{\Omega_x}|\nabla_{\mathbb{R}^2}(u\circ\psi_x)|^2dy=\int_{U_x}|\nabla_gu|^2dv_g\leq \lambda_{\rm N}(\Sigma)/
(\lambda_{\rm N}(\Sigma)-\alpha)\ee
and that
\be\label{int}\int_{{\Omega_x}}(u\circ\psi_x)^2dy\leq \le(\max_{\overline{\Omega_x}}\exp(2f)\ri)\int_{U_x}u^2dv_g\leq
\f{\max_{\overline{\Omega_x}}\exp(2f)}{\lambda_{\rm N}(\Sigma)-\alpha}.\ee
Combining (\ref{euc}), (\ref{int}) and Chang-Yang's result (\ref{C-Y}), we conclude that there must be two constants $\gamma_x<2\pi$ and $C_x>0$
satisfying
$$\int_{\Omega_x}\exp(\gamma_x(u\circ\psi_x)^2)dy\leq C_x.$$
It then follows that
\be\label{U}\int_{U_x}\exp(\gamma_xu^2)dv_g\leq C_x\max_{\overline{\Omega_x}}\exp(-2f).\ee
In the case $x\in\p\Sigma$, the estimate (\ref{U}) still holds for some constant $C_x$ due to Lemma \ref{prop} and Chang-Yang's result (\ref{C-Y}).
Since $(\overline{\Sigma},g)$ is compact, we can choose $\ell$ sets $\{U_{x_i}\}_{i=1}^\ell$ satisfying $\cup_{i=1}^\ell \overline{U_{x_i}}\supset
\overline{\Sigma}$,
  where $U_{x_i}$ is given as above. We immediately get the desired result. $\hfill\Box$\\

In view of Lemma \ref{lemma1}, we let
\be\label{gamma}\gamma^\ast=\sup\le\{\gamma: \sup_{u\in\mathscr{H},\,\|u\|_{1,\alpha}\leq 1}\int_\Sigma \exp(\gamma u^2)dv_g<+\infty\ri\}.\ee

\begin{lemma}\label{lemma2}
There holds $\gamma^\ast\leq 2\pi$.
\end{lemma}

{\it Proof}.
Recall the Moser function sequence \cite{Moser}
\be\label{Moser-sequence}M_{k}(y,r)=\le\{
\begin{array}{lll}
\sqrt{\f{\log k}{4\pi}}&{\rm when}&|y|\leq rk^{-1/4}\\[1.5ex]
\sqrt{\f{4}{\pi\log k}}\log\f{r}{|y|}&{\rm when}&rk^{-1/4}<|y|\leq r\\[1.5ex]
0&{\rm when}& |y|>r
\end{array}
\ri.
\ee
for all $y\in\mathbb{R}^2$, $r>0$ and $k\in\mathbb{N}$. It can be checked that
\be\label{te-1}\int_{\mathbb{B}_r^+}|\nabla_{\mathbb{R}^2}M_k(y,r)|^2dy=1,\ee
\be\label{te-2}\int_{\mathbb{B}_r^+}M_k(y,r)dy=o_k(1)+o_r(1),\ee
that
\be\label{te-5}\int_{\mathbb{B}_r^+}M_k^2(y,r)dy=o_k(1)+o_r(1),\ee
where $o_k(1)\ra 0$ as $k\ra\infty$, $o_r(1)\ra 0$ as $r\ra0$, and that
\be\label{te-3}\int_{\mathbb{B}_r^+}\exp(\gamma M_k^2(y,r))dy\geq \int_{\mathbb{B}_{rk^{-1/4}}^+}\exp(\gamma M_k^2(y,r))dy
=\f{\pi}{2}r^2k^{\f{\gamma}{4\pi}-\f{1}{2}}.\ee

Now we fix a point $p\in\p\Sigma$ and choose an isothermal coordinate system $(\overline{U_p},\psi;\{y_1,y_2\})$
near $p$, where $\overline{U_p}=\psi^{-1}(\overline{\mathbb{B}_r^+})\subset \overline{\Sigma}$ for some $r>0$.
In this coordinate system, the metric
$$g=\exp(2f(y))(dy_1^2+dy_2^2),$$
where $f\in C^1(\overline{\mathbb{B}_r^+})$ with $f(0,0)=0$.
Define a sequence of functions
$$\widetilde{M}_k(x,r)=\le\{
\begin{array}{lll} M_k(\psi(x),r)&{\rm if}&x\in \overline{U_p}\\[1.5ex]
0&{\rm if}& x\in\overline{\Sigma}\setminus \overline{U_p}.
\end{array}
\ri.$$
In view of (\ref{Moser-sequence})-(\ref{te-5}), we have that
$$\int_\Sigma|\nabla_g\widetilde{M}_k|^2dv_g=1,\,\int_\Sigma \widetilde{M}_kdv_g=o_k(1)+o_r(1),\,\int_\Sigma \widetilde{M}_k^2dv_g=o_k(1)+o_r(1).$$
Let
$$Q_k=\widetilde{M}_k-\f{1}{{\rm Area}(\Sigma)}\int_\Sigma\widetilde{M}_kdv_g.$$
It follows that $Q_k\in\mathscr{H}$ and $\|Q_k\|_{1,\alpha}=1+o_k(1)+o_r(1)$. This together with  (\ref{te-3}) implies
\bna
\int_\Sigma\exp(\gamma Q_k^2/\|Q_k\|_{1,\alpha}^2)dv_g\geq (1+o_r(1))\f{\pi}{2}r^2k^{\f{\gamma}{4\pi}-\f{1}{2}+o_k(1)+o_r(1)}.
\ena
Therefore if $\gamma>2\pi$, then we have by choosing sufficiently small $r$ and passing to the limit $k\ra\infty$,
$$\int_\Sigma\exp(\gamma Q_k^2/\|Q_k\|_{1,\alpha}^2)dv_g\ra+\infty.$$
This leads to $\gamma^\ast\leq 2\pi$. $\hfill\Box$\\

Furthermore, we have

\begin{lemma}\label{lemma3} There holds
$\gamma^\ast= 2\pi$.
\end{lemma}
{\it Proof}.
In view of Lemma \ref{lemma2}, we only need to show $\gamma^\ast$ can not be strictly less than $2\pi$.
By the definition of $\gamma^\ast$ (see (\ref{gamma}) above), there exists a function sequence $(w_j)\subset\mathscr{H}$ with $\|w_j\|_{1,\alpha}\leq 1$
such that
\be\label{infty}\int_\Sigma \exp((\gamma^\ast+{1}/{j})w_j^2)dv_g\ra +\infty\ee
as $j\ra\infty$. Clearly, there exists some $w\in\mathscr{H}$ with $\|w\|_{1,\alpha}\leq1$ such that $w_j\rightharpoonup w$ weakly in $W^{1,2}(\Sigma,g)$, $w_j\ra w$ strongly in $L^q(\Sigma)$ for any $q>0$
and $w_j\ra w$ almost everywhere in $\Sigma$. We now {\it claim} $w\equiv 0$ in $\Sigma$. Supposing the contrary, we would have
$$\|w_j-w\|_{1,\alpha}^2=1-\|w\|_{1,\alpha}^2+o_j(1)\leq 1-\f{1}{2}\|w\|_{1,\alpha}^2$$
for sufficiently large $j$. For any $\epsilon>0$, one has by using the Young inequality,
$ab\leq \epsilon a^2+b^2/(4\epsilon)$,  and the H\"older inequality that
$$\int_\Sigma\exp((\gamma^\ast+{1}/{j})w_j^2)dv_g\leq C\le(\int_\Sigma\exp\le((\gamma^\ast+1/j)(1+2\epsilon)
\|w_j-w\|_{1,\alpha}^2\f{(w_j-w)^2}{\|w_j-w\|_{1,\alpha}^2}\ri)dv_g\ri)^{\f{1+\epsilon}{1+2\epsilon}}$$
for some constant $C$ depending only on $\epsilon$ and $w$. Taking $\epsilon$ such that
$1+2\epsilon=(1-\|w\|_{1,\alpha}^2/3)/(1-\|w\|_{1,\alpha}^2/2)$, we have
$$(\gamma^\ast+1/j)(1+2\epsilon)\|w_j-w\|_{1,\alpha}^2\leq (1-\|w\|_{1,\alpha}^2/3)(\gamma^\ast+1/j)\leq
(1-\|w\|_{1,\alpha}^2/4)\gamma^\ast,$$
provided that $j\geq j_0$ for sufficiently large $j_0$. As a consequence
$$\int_\Sigma\exp((\gamma^\ast+{1}/{j})w_j^2)dv_g\leq \int_\Sigma \exp\le((1-\|w\|_{1,\alpha}^2/4)\gamma^\ast
\f{(w_j-w)^2}{\|w_j-w\|_{1,\alpha}^2}\ri)dv_g\leq C$$
for some constant $C$. This contradicts (\ref{infty}) and confirms our claim $w\equiv 0$.

Suppose that $\gamma^\ast<2\pi$. Similarly as in the proof of Lemma \ref{lemma1}, for any $x\in\overline{\Sigma}$,
we choose an isothermal coordinate system $(U_x,\psi_x;\{y_1,y_2\})$,
where $\psi_x: U_x\ra\Omega_x\subset\mathbb{R}^2$ is a diffeomorphism. In such a coordinate system,
the metric $g=\exp(2f(y))(dy_1^2+dy_2^2)$, where $f\in C^1(\overline{\Omega_x})$ with
$f(0,0)=0$. By the above consideration, $w_j$ converges to $0$ strongly in $L^q(\Sigma)$ for any $q>0$. It follows that
$$\int_{\Omega_x}w_j\circ \psi_x^{-1}dy=o_j(1)$$
and
$$\int_{\Omega_x}|\nabla_{\mathbb{R}^2}(w_j\circ \psi_x^{-1})|^2dy\leq 1+o_j(1).$$
Hence for any $\gamma$, $\gamma^\ast<\gamma<2\pi$, we have by Chang-Yang's result (\ref{C-Y}),
$$\int_{\Omega_x}\exp(\gamma (w_j\circ\psi_x)^2)dy\leq C.$$
Similarly we have
$$\int_{U_x}\exp(\gamma w_j^2)dv_g\leq C.$$
Since $\overline{\Sigma}$ is compact, by choosing finitely many isothermal coordinate systems covering $\overline{\Sigma}$, we
conclude
$$\int_\Sigma\exp(\gamma w_j^2)dv_g\leq C$$
for some constant $C$ depending only on $(\Sigma,g)$, $\gamma^\ast$ and $\gamma$.
This contradicts (\ref{infty}) and concludes that $\gamma^\ast$ must be $2\pi$. $\hfill\Box$

  \subsection{Existence of extremals for subcritical Trudinger-Moser functionals}\label{4.2}

  Using a direct method of variation, we can prove the attainability of the supremum in (\ref{ineq-1})
  in the case $\gamma<2\pi$.
  In particular, we have  the following:

  \begin{lemma}\label{subcritical}
  For any $k\in\mathbb{N}$, there exists a $u_k\in\mathscr{H}\cap C^1(\overline{\Sigma})$ with $\|u_k\|_{1,\alpha}=1$ such that
  $$\int_\Sigma \exp((2\pi-1/k)u_k^2)dv_g=\sup_{u\in\mathscr{H},\,\|u\|_{1,\alpha}\leq 1}\int_\Sigma \exp((2\pi-1/k)u^2)dv_g.$$
  Moreover $u_k$ satisfies the following Euler-Lagrange equation
  \be\label{E-L}\le\{\begin{array}{lll}
  \Delta_gu_k-\alpha u_k=\f{1}{\lambda_k}u_k\exp(\gamma_ku_k^2)-\f{\mu_k}{\lambda_k}&{\rm in}&\Sigma\\[1.5ex]
  {\p u_k}/{\p {\bm\nu}}=0&{\rm on}&\p\Sigma\\[1.5ex]
  \lambda_k=\int_\Sigma u_k^2\exp(\gamma_ku_k^2)dv_g\\[1.5ex]
  \mu_k=\f{1}{{\rm Area}(\Sigma)}\int_\Sigma u_k\exp(\gamma_ku_k^2)dv_g\\[1.5ex]
  \gamma_k=2\pi-1/k,
  \end{array}\ri.\ee
  where $\bm\nu$ denotes the unit outward vector field on $\p\Sigma$.
  \end{lemma}
  {\it Proof}.
  The proof is based on a direct variational method. Take a function sequence $(u_j)\subset\mathscr{H}$ satisfying $\|u_j\|_{1,\alpha}\leq 1$
  and
  \be\label{lim}\lim_{j\ra\infty}\int_\Sigma \exp((2\pi-1/k)u_j^2)dv_g=\sup_{u\in\mathscr{H},\,\|u\|_{1,\alpha}\leq 1}\int_\Sigma \exp((2\pi-1/k)u^2)dv_g.\ee
  Since $\alpha<\lambda_{\rm N}(\Sigma)$, $(u_j)\subset\mathscr{H}$ is bounded in $W^{1,2}(\Sigma)$. Hence one can find $u_k\in\mathscr{H}$ such
  that $u_j$ converges to $u_k$ weakly in $W^{1,2}(\Sigma,g)$, strongly in $L^q(\Sigma)$ for all $q>0$, and almost everywhere in $\Sigma$.
  It then follows that $\|u_k\|_{1,\alpha}\leq 1$. By Lemma \ref{lemma3}, $\exp((2\pi-1/k)u_j^2)$ is bounded in $L^s(\Sigma)$ for some $s>1$,
  and thus $\exp((2\pi-1/k)u_j^2)$ converges to $\exp((2\pi-1/k)u_k^2)$ in $L^1(\Sigma)$. This together with (\ref{lim}) leads to
  \be\label{sup}\int_\Sigma \exp((2\pi-1/k)u_k^2)dv_g=\sup_{u\in\mathscr{H},\,\|u\|_{1,\alpha}\leq 1}\int_\Sigma \exp((2\pi-1/k)u^2)dv_g.\ee
  Now we show $\|u_k\|_{1,\alpha}=1$. Suppose $\|u_k\|_{1,\alpha}<1$. Then
  $$\int_\Sigma \exp((2\pi-1/k)u_k^2/\|u_k\|_{1,\alpha}^2)dv_g>\int_\Sigma \exp((2\pi-1/k)u_k^2)dv_g.$$
  This contradicts (\ref{sup}) and implies that $\|u_k\|_{1,\alpha}=1$.
  A simple calculation shows the Euler-Lagrange equation of $u_k$ is (\ref{E-L}).
  Applying elliptic estimates to (\ref{E-L}), we have $u_k\in C^1(\overline{\Sigma})$.
  $\hfill\Box$

  \subsection{Blow-up analysis on the boundary}\label{4.3}

  Let $c_k=\max_{\overline{\Sigma}}|u_k|$. If $c_k$ is bounded, then applying elliptic estimates to (\ref{E-L}), we conclude that there exists some
  $u^\ast\in\mathscr{H}$ with $\|u^\ast\|_{1,\alpha}=1$ such that $u_k\ra u^\ast$ in $C^1(\overline{\Sigma})$. Clearly $u^\ast$ is the desired extremal
  function. In the following, noting that $-u_k$ is also a solution of (\ref{E-L}), we assume without loss of generality that
  \be\label{c-} c_k=u_k(x_k)=\max_{\overline{\Sigma}}|u_k|\ra+\infty\ee
  and
  \be\label{lim-point}x_k\ra x_0\in\overline{\Sigma}\ee
  as $k\ra\infty$. Since $u_k$ is bounded in $W^{1,2}(\Sigma,g)$, we assume up to a subsequence
  $u_k$ converges to $u_0\in\mathscr{H}$ weakly in $W^{1,2}(\Sigma,g)$, strongly in $L^q(\Sigma,g)$ for all $q>0$,
  and almost everywhere in $\Sigma$.
  \begin{lemma}\label{concern}
   $u_0\equiv 0$, $x_0\in\p\Sigma$, and $|\nabla_gu_k|^2dv_g\rightharpoonup \delta_{x_0}$ in the sense of measure.
  \end{lemma}
  {\it Proof}.
  Firstly we prove $u_0\not\equiv 0$. Suppose not. There holds
  $$\|u_k-u_0\|_{1,\alpha}^2=1-\|u_0\|_{1,\alpha}^2+o_k(1)\leq 1-\f{1}{2}\|u_0\|_{1,\alpha}^2$$
  for sufficiently large $k$. By the inequality $(a+b)^2\leq (1+\epsilon)a^2+(1+1/\epsilon)b^2$, the H\"older inequality,
  and Lemma \ref{lemma3},
  we have that $\exp(\gamma_ku_k^2)$ is bounded in $L^q(\Sigma,g)$ for some $q>1$. Then applying elliptic estimates to (\ref{E-L}),
  we conclude that $u_k$ is uniformly bounded in $\overline{\Sigma}$, which contradicts (\ref{c-}). Hence $u_0\equiv 0$.

  Secondly, in view of (\ref{lim-point}), we show $x_0\in\p\Sigma$. Suppose $x_0$ is an inner point of $\Sigma$. Choose an isothermal coordinate system
  $(U_{x_0},\psi_{x_0};\{y_1,y_2\})$ around $x_0$ such that $\psi_{x_0}(U_{x_0})=\mathbb{B}_{r_0}\subset\mathbb{R}^2$.
  In this coordinate system, the metric $g$ can be written as $g=\exp(2f(y))(dy_1^2+dy_2^2)$, where $f\in C^1(\overline{\mathbb{B}_{r_0}})$
  with $f(0,0)=0$.
  Take a cut-off function $\phi\in C_0^2(\mathbb{B}_{r_0})$ satisfying $0\leq\phi\leq 1$ and $\phi\equiv 1$ on $\mathbb{B}_{r_0/2}$.
  One has $\phi(u_k\circ\psi_{x_0}^{-1})\in W_0^{1,2}(\mathbb{B}_{r_0})$ and $\|\nabla_{\mathbb{R}^2}(\phi(u_k\circ\psi_{x_0}^{-1}))\|_2^2\leq
  1+o_k(1)$. Thus Moser's inequality (\ref{Moser}) leads to
  $$\int_{\mathbb{B}_{r_0/2}}\exp\le(2\pi q (u_k\circ\psi_{x_0}^{-1})^2\ri)dx\leq\int_{\mathbb{B}_{r_0}}\exp\le(2\pi q \phi^2(u_k\circ\psi_{x_0}^{-1})^2\ri)dx\leq C$$
  for some $q>1$ and constant $C$. This immediately implies that
  ${\lambda_k^{-1}}(u_k\exp(\gamma_ku_k^2)-\lambda_k^{-1}{\mu_k}$
  is bounded in $L^{q^\prime}(\psi_{x_0}^{-1}(\mathbb{B}_{r_0/2}))$ for some $1<q^\prime<q$. Applying elliptic estimates to (\ref{E-L}), we conclude
  that $u_k$ is uniformly bounded in  $\psi_{x_0}^{-1}(\mathbb{B}_0{(r_0/4)})$, contradicting (\ref{c-}). Therefore $x_0\in\p\Sigma$.

  As for the final assertion, we first {\it claim} the following
  \be\label{one}\lim_{r\ra 0}\lim_{k\ra\infty}\int_{B_{r}(x_0)}|\nabla_gu_k|^2dv_g=1,\ee
  where $B_{r}(x_0)\subset\Sigma$ denotes the geodesic ball centered at $x_0$ with radius $r$.
  For otherwise, there exist $a<1$, $r>0$ and $k_0>0$ such that
  $$\int_{B_{r}(x_0)}|\nabla_gu_k|^2dv_g\leq a,\quad\forall k\geq k_0.$$
  Then similarly as we derive $x_0\in\p\Sigma$, we conclude that $u_k$ is uniformly bounded in $B_{r/2}(x_0)$, which
  contradicts (\ref{c-}) again. Hence (\ref{one}) holds. For any $\varphi\in C^0(\overline{\Sigma})$, we have
  \bna\lim_{k\ra\infty}\int_{\Sigma}\varphi|\nabla_gu_k|^2dv_g&=&\lim_{k\ra\infty}\le(\int_{\Sigma\setminus B_{r}(x_0)}
  \varphi|\nabla_gu_k|^2dv_g+\int_{B_{r}(x_0)}
  \varphi|\nabla_gu_k|^2dv_g\ri)\\
  &=&\varphi(x_0),\ena
  which is the desired result. $\hfill\Box$\\

  Let $x_0$ be given as in Lemma \ref{concern}. From now on until the end of this section, we use the isothermal coordinate system
  \be\label{iso-coord}(\overline{U_{x_0}},\psi_{x_0};\{y_1,y_2\}),\,\,\psi_{x_0} (U_{x_0})=\mathbb{B}_{r_0}^+,
  \,\,\psi_{x_0}(\overline{U_{x_0}}\cap\p\Sigma)=\p\mathbb{B}_{r_0}^+\cap\p\mathbb{R}^{2+},\ee
  \be\label{metr-repre}\psi_{x_0}(x_0)=(0);\,\,g=\exp(2f(y))(dy_1^2+dy_2^2),\,f(0,0)=0,\,f\in C^1(\overline{\mathbb{B}_{r_0}^+});\ee
  moreover, in this coordinate system,
   the unit outward vector field $\bm\nu$ on the boundary
  $\p\Sigma$ can be written as
  $\bm\nu=\exp(-f(y)){\p}/{\p y_2}$. For any $u\in C^1(\overline{\Sigma})$, the normal derivative $\p u/\p\bm{\nu}$ can be represented by
  \be\label{norm}\f{\p u}{\p\bm\nu}=\exp(-f(y))\f{\p}{\p y_2}(u\circ\psi_{x_0}^{-1}).\ee
  For simplicity we write
  $$f_k={\lambda_k^{-1}}u_k\exp(\gamma_ku_k^2)-{\lambda_k^{-1}}{\mu_k}+\alpha u_k,$$
  where $\alpha$, $u_k$, $\mu_k$, $\gamma_k$ and $\lambda_k$ are defined as in (\ref{E-L}).
  We set
  \be\label{u-ast}u_k^\ast(y)=\le\{\begin{array}{lll}
  u_k\circ\psi_{x_0}^{-1}(y)&{\rm if}&y\in\overline{\mathbb{B}_{r_0}^+}\\[1.5ex]
  u_k\circ\psi_{x_0}^{-1}(y_1,-y_2)&{\rm if}&y\in\mathbb{B}_{r_0}^-\end{array}
  \ri.\ee
  and
  $$f_k^\ast(y)=\le\{\begin{array}{lll}
  \exp(2f(y))(f_k\circ\psi_{x_0}^{-1})(y)&{\rm if}&y\in\overline{\mathbb{B}_{r_0}^+}\\[1.5ex]
  \exp(2f(y_1,-y_2))(f_k\circ\psi_{x_0}^{-1})(y_1,-y_2)&{\rm if}&y\in\mathbb{B}_{r_0}^-,\end{array}
  \ri.$$
  where $\mathbb{B}_{r_0}^-=\mathbb{B}_{r_0}\setminus\overline{\mathbb{B}_{r_0}^+}$.
  Define two function sequences $\phi_k:\mathbb{B}_R\ra \mathbb{R}$ and $\eta_k:\mathbb{B}_R\ra \mathbb{R}$
  for any fixed $R>0$ by
  \be\label{scaling}\phi_k(z)=\f{u_k^\ast(\widetilde{x}_k+r_kz)}{c_k},\quad
  \eta_k(z)=c_k(u_k^\ast(\widetilde{x}_k+r_kz)-c_k),\ee
  where $\widetilde{x}_k=\psi_{x_0}(x_k)$, $c_k$ is defined as in (\ref{c-}) and $r_k>0$ satisfies
  \be\label{rk}r_k^2=\f{\lambda_k}{c_k^2}\exp(-\gamma_k c_k^2).\ee

  \begin{lemma}\label{bubble}
  $\phi_k$ and $\eta_k$ are distributional solutions of
  \be\label{p-equ}-\Delta_{\mathbb{R}^2}\phi_k(z)=\f{r_k^2}{c_k}f_k^\ast(\widetilde{x}_k+r_kz)\quad{\rm in}\quad \mathbb{B}_R\ee
  and
  \be\label{e-equ}-\Delta_{\mathbb{R}^2}\eta_k(z)=r_k^2c_kf_k^\ast(\widetilde{x}_k+r_kz)\quad{\rm in}\quad \mathbb{B}_R\ee
  respectively.
  \end{lemma}
  {\it Proof}. In view of (\ref{u-ast}), $u_k^\ast\in W^{1,2}(\mathbb{B}_{r_0})\cap C^1(\overline{\mathbb{B}_{r_0}^+})$. Since $\p u_k/\p\bm{\nu}=0$ on $\p\Sigma$, in view of (\ref{norm}), we have $\p u_k^\ast/\p y_2=0$
  on $\p\mathbb{B}_{r_0}^+\cap\p\mathbb{R}^{2+}$. We {\it claim} that $u_k^\ast$ is a distributional solution of
  the equation
  \be\label{symm}-\Delta_{\mathbb{R}^2}u_k^\ast=f_k^\ast\quad{\rm in}\quad\mathbb{B}_{r_0}.\ee
  To see this, for any $\varphi\in C_0^\infty(\mathbb{B}_{r_0})$, we obtain
  \bna
  -\int_{\mathbb{B}_{r_0}}u_k^\ast\Delta_{\mathbb{R}^2}\varphi dy&=&\int_{\mathbb{B}_{r_0}}\nabla_{\mathbb{R}^2}u_k^\ast\nabla_{\mathbb{R}^2}\varphi dy\\&=&\int_{\mathbb{B}_{r_0}^+}\nabla_{\mathbb{R}^2}u_k^\ast\nabla_{\mathbb{R}^2}\varphi dy+
  \int_{\mathbb{B}_{r_0}\setminus\mathbb{B}_{r_0}^+}\nabla_{\mathbb{R}^2}u_k^\ast\nabla_{\mathbb{R}^2}\varphi dy\\
  &=&\int_{\p\mathbb{B}_{r_0}^+\cap\p\mathbb{R}^{2+}}\f{\p u_k^\ast}{\p y_2}\varphi dy-\int_{\mathbb{B}_{r_0}^+} (\Delta_{\mathbb{R}^{2}}u_k^\ast)\varphi dy\\&&+\int_{\p\mathbb{B}_{r_0}^+\cap\p\mathbb{R}^{2+}}\f{\p u_k^\ast}{\p y_2}\varphi dy-\int_{\mathbb{B}_{r_0}\setminus\mathbb{B}_{r_0}^+} (\Delta_{\mathbb{R}^{2}}u_k^\ast)\varphi dy\\
  &=&\int_{\mathbb{B}_{r_0}}f_k^\ast\varphi dy,
    \ena
    which concludes that $u_k^\ast$ satisfies (\ref{symm}) in the distributional sense.

  We next prove that $\phi_k$ is a distributional solution of (\ref{p-equ}). Let $R>0$ be fixed. For any $\psi\in C_0^\infty(\mathbb{B}_R)$,
  we denote $\widetilde{\psi}(y)=\psi(\widetilde{x}_k+r_kz)$. Obviously $\widetilde{\psi}\in C_0^\infty(\mathbb{B}_{Rr_k}(\widetilde{x}_k))$.
  Since $u_k^\ast$ is a distributional solution of (\ref{symm}), it then follows that
  \bna
  \int_{\mathbb{B}_R}\phi_k(z) \Delta_{\mathbb{R}^2}\psi(z)dz&=&\int_{\mathbb{B}_{Rr_k}(\widetilde{x}_k)}\f{u_k^\ast(y)}{c_k}\f{1}{r_k^2}
  \Delta_{\mathbb{R}^2}\psi\le(\f{y-\widetilde{x}_k}{r_k}\ri)dy\\&=&\int_{\mathbb{B}_{Rr_k}(\widetilde{x}_k)}\f{u_k^\ast(y)}{c_k}
  \Delta_{\mathbb{R}^2}\widetilde{\psi}(y)dy\\
  &=&\int_{\mathbb{B}_{Rr_k}(\widetilde{x}_k)}\f{f_k^\ast(y)}{c_k}\widetilde{\psi}(y)dy\\
  &=&\int_{\mathbb{B}_R}\f{f_k^\ast(\widetilde{x}_k+r_kz)}{c_k}\psi(z)r_k^2dz.
  \ena
  Hence $\phi_k$ satisfies (\ref{p-equ}) in the distributional sense.

  In the same way, it can be proved that $\eta_k$ is a distributional solution of (\ref{e-equ}).  $\hfill\Box$

  \begin{lemma}
  For any $\nu<2\pi$, there holds $r_k\exp(\nu c_k^2)$ converges to $0$ as $k\ra\infty$.
  \end{lemma}
  {\it Proof}. Using the H\"older inequality, the fact $u_k\ra 0$ strongly in $L^q(\Sigma,g)$ for any $q>0$, and Lemma \ref{lemma3},
  we have for any $\nu<2\pi$,
  $$\int_\Sigma u_k^2\exp(\nu u_k^2)dv_g=o_k(1).$$
  This together with (\ref{rk}) and the definition of $\lambda_k$ (see (\ref{E-L}) above) gives the desired result. $\hfill\Box$

  \begin{lemma}\label{lemma9} Let $\phi_k$ and $\eta_k$ be defined as in (\ref{scaling}). Then
  $\phi_k\ra\phi_0$ and $\eta_k\ra\eta_0$ in $C^1_{\rm loc}(\mathbb{R}^2)$, where
  $\phi_0(z)\equiv 1$ and
  $$\eta_0(z)=-\f{1}{2\pi}\log(1+\f{\pi}{2}|z|^2),\quad\forall z\in\mathbb{R}^2.$$
  \end{lemma}

  {\it Proof}. The proof is based on the elliptic estimates on (\ref{p-equ}) and (\ref{e-equ}). We omit the details but refer the reader to \cite{LiJPDE,Yang-Tran}.$\hfill\Box$

  \begin{lemma}\label{lemma10}
  Let $\widetilde{x}_k=(y_{1,k},y_{2,k})$. Then $y_{2,k}/r_k\ra 0$ as $k\ra\infty$.
  \end{lemma}
  {\it Proof}. Without loss of generality, we assume as $k\ra\infty$,
  \be\label{limt}y_{2,k}/r_k\ra \ell\ee
  for some $\ell>0$. Noting that under the change of variable $y=\widetilde{x}_k+r_kz$, the set $\mathbb{B}_{Rr_k}(\widetilde{x}_k)
  \cap\mathbb{R}^{2+}=\{y=(y_1,y_2)\in
  \mathbb{B}_{Rr_k}(\widetilde{x}_k): y_2>0\}$ is mapped onto $B_{\ell,k}=\{z=(z_1,z_2)\in
  \mathbb{B}_R: z_2>-y_{2,k}/r_k=-\ell(1+o_k(1))\}$, we calculate by noticing (\ref{limt}),
  \bea\nonumber
  1&=&\int_\Sigma \f{1}{\lambda_k}u_k^2\exp(\gamma_ku_k^2)dv_g\\\nonumber
  &\geq&\int_{\mathbb{B}_{Rr_k}(\widetilde{x}_k)
  \cap\mathbb{R}^{2+}}\f{1}{\lambda_k}{u_k^\ast}^2\exp(\gamma_k{u_k^\ast}^2)\exp(2v(y))dy\\\nonumber
  &=&(1+o_k(1))\int_{B_{\ell,k}}\phi_k^2(z)\exp(\gamma_k(1+\phi_k(z))\eta_k(z))dz
  \\\label{00}
  &=&(1+o_k(1))\int_{\mathbb{B}_R\cap\{z_2>-\ell\}}\exp(4\pi\eta_0(z))dz.
  \eea
  Note that
  $$\int_{z_2>-\ell}\exp(4\pi \eta_0)dz>\int_{\mathbb{R}^{2+}}\exp(4\pi \eta_0)dz=1.$$
  By passing to the limit $k\ra\infty$ and then $R\ra\infty$ in (\ref{00}), we get a contradiction.
  This ends the proof of the lemma. $\hfill\Box$\\

  Let $\widetilde{x}_{0,k}=(y_{1,k},0)$. We define two function sequences modified from (\ref{scaling}) by
  \be\label{x-0k}\phi_{1,k}(z)=\f{u_k^\ast(\widetilde{x}_{0,k}+r_kz)}{c_k},\quad \eta_{1,k}(z)=
  c_k(u_k^\ast(\widetilde{x}_{0,k}+r_kz)-c_k)\ee
  for $z\in\Omega_k=\{z\in\mathbb{R}^2: \widetilde{x}_{0,k}+r_kz\in \mathbb{B}_{r_0}\}$.

  \begin{lemma}\label{x'}
  $\phi_{1,k}\ra\phi_0$ and $\eta_{1,k}\ra\eta_0$  in $C^1_{\rm loc}(\mathbb{R}^2)$, where $\phi_0$ and $\eta_0$ are given as in Lemma \ref{lemma9}.
  \end{lemma}

  {\it Proof}. Note that
  \bna\widetilde{x}_{0,k}+r_kz&=&\widetilde{x}_k+r_k\le(z+\f{\widetilde{x}_{0,k}-\widetilde{x}_{k}}{r_k}\ri)\\
  &=&\widetilde{x}_k+r_k\le(z+(0,{y_{2,k}}/{r_k})\ri).\ena
  Then the lemma follows from Lemmas \ref{lemma9} and \ref{lemma10}. $\hfill\Box$\\

  For any $0<\beta<1$,  let $u_{k,\beta}=\min\{u_k,\beta c_k\}$.  Similar to \cite{LiJPDE},  we shall show
  \be\label{uk-beta}\lim_{k\ra\infty}\int_\Sigma |\nabla_gu_{k,\beta}|^2dv_g=\beta.\ee
  To this end, since $\p u_k/\p{\bm\nu}=0$ on $\p\Sigma$, we have by the divergence theorem, Lemmas \ref{lemma9}
  and \ref{lemma10},
  \bna\nonumber
  \int_\Sigma|\nabla_gu_{k,\beta}|^2dv_g&=&\int_\Sigma u_{k,\beta}\Delta_gu_kdv_g\\\nonumber
  &=&\int_\Sigma u_{k,\beta}\f{1}{\lambda_k}u_k\exp(\gamma_ku_k^2)dv_g+o_k(1)\\\nonumber
  &\geq&(1+o_k(1))\int_{\mathbb{B}_{Rr_k}(\widetilde{x}_k)\cap \mathbb{B}_{r_0}^+} \f{\beta}{\lambda_k}c_ku_k^\ast\exp(\gamma_k{u_k^\ast}^2)
  dy+o_k(1)\\\nonumber
  &=&(\beta+o_k(1))\int_{\mathbb{B}_R\cap \{z_2>-y_{2,k}/r_k\}}\exp(4\pi \eta_0(z))dz+o_k(1)\\
  &=&\beta\int_{\mathbb{B}_R\cap\mathbb{R}^{2+}}\exp(4\pi \eta_0(z))dz+o_k(1).
  \ena
  Letting $k\ra\infty$ first, and then $R\ra\infty$, we have
  \be\label{geq}\liminf_{k\ra\infty}\int_\Sigma|\nabla_gu_{k,\beta}|^2dv_g\geq \beta.\ee
  In the same way, we estimate
  \bna
  \int_\Sigma |\nabla_g (u_k-u_{k,\beta})^+|^2dv_g&=&\int_\Sigma (u_k-u_{k,\beta})^+\Delta_g u_kdv_g\\
  &\geq&(1-\beta)\int_{\mathbb{B}_R\cap\mathbb{R}^{2+}}\exp(4\pi \eta_0(z))dz+o_k(1).
  \ena
  Then we get an analog of (\ref{geq}), namely
   \be\label{geq-1}\liminf_{k\ra\infty}\int_\Sigma|\nabla_g(u_k-u_{k,\beta})^+|^2dv_g\geq 1-\beta.\ee
  Hence the equality
  $$\int_\Sigma |\nabla_gu_{k,\beta}|^2dv_g+\int_\Sigma |\nabla_g(u_k-u_{k,\beta})^+|^2dv_g=\int_\Sigma |\nabla_gu_{k}|^2dv_g=
  1+\alpha\int_\Sigma u_k^2dv_g$$
  together with (\ref{geq}), (\ref{geq-1}) and Lemma \ref{concern} implies (\ref{uk-beta}).

  \begin{lemma}\label{jixian} Under the assumption $c_k\ra\infty$, there holds
  \be\label{L1}\sup_{u\in\mathscr{H},\,\|u\|_{1,\alpha}\leq 1}\int_\Sigma\exp(2\pi u^2)dv_g={\rm Area}(\Sigma)+\lim_{k\ra\infty}\f{\lambda_k}{c_k^2}.\ee
  As a consequence,
  \be\label{L0}c_k/\lambda_k\ra 0\quad {\rm as}\quad k\ra\infty.\ee
  \end{lemma}
  {\it Proof}. Note that
  \be\label{sup-2pi}\lim_{k\ra\infty}\int_\Sigma \exp(\gamma_ku_k^2)dv_g=\sup_{u\in\mathscr{H},\,\|u\|_{1,\alpha}\leq 1}\int_\Sigma\exp(2\pi u^2)dv_g.\ee
  Given any $0<\beta<1$. On one hand, we have by (\ref{uk-beta}),
  \bea\int_\Sigma \exp(\gamma_ku_k^2)dv_g&\geq&\int_{u_k> \beta c_k}\f{u_k^2}{c_k^2}\exp(\gamma_ku_k^2)dv_g+
  \int_{u_k\leq \beta c_k}\exp(\gamma_ku_k^2)dv_g\nonumber\\
  &=&\f{\lambda_k}{c_k^2}-\int_{u_k\leq \beta c_k}\f{u_k^2}{c_k^2}\exp(\gamma_ku_k^2)dv_g+
  \int_{u_k\leq \beta c_k}\exp(\gamma_ku_k^2)dv_g\nonumber \\
  &=&\f{\lambda_k}{c_k^2}+{\rm Area}(\Sigma)+o_k(1).\label{sup-2}\eea
  On the other hand, we also obtain by using (\ref{uk-beta}),
  \bea\nonumber\int_\Sigma \exp(\gamma_ku_k^2)dv_g&\leq&\int_{u_k> \beta c_k}\f{u_k^2}{\beta^2c_k^2}\exp(\gamma_ku_k^2)dv_g+
  \int_{u_k\leq \beta c_k}\exp(\gamma_ku_k^2)dv_g\\\label{sup-3}
  &\leq&\f{1}{\beta^2}\f{\lambda_k}{c_k^2}+{\rm Area}(\Sigma)+o_k(1).\eea
  Combining (\ref{sup-2pi})-(\ref{sup-3}), we get (\ref{L1}) by passing to the limit $k\ra\infty$ first,
  and then $\beta\ra 1$.

  For the second assertion, we suppose the contrary, there exists some constant $\varrho>0$ such that up to a subsequence,
  $c_k/\lambda_k\geq\varrho$. Hence $\lambda_k/c_k^2\leq 1/(\varrho c_k)=o_k(1)$, which together with (\ref{L1}) leads to
  $$\sup_{u\in\mathscr{H},\,\|u\|_{1,\alpha}\leq 1}\int_\Sigma\exp(2\pi u^2)dv_g={\rm Area}(\Sigma),$$
  which is impossible. Therefore (\ref{L0}) holds. $\hfill\Box$

  \begin{lemma}\label{tends-delta}
  For any $\varphi\in C^2(\overline{\Sigma})$, there holds
  \be\label{det}\int_\Sigma\varphi\f{c_k }{\lambda_k}|u_k|\exp(\gamma_ku_k^2)dv_g=\varphi(x_0)+o_k(1),\,\,\,
  \int_\Sigma\varphi\f{c_k }{\lambda_k}u_k\exp(\gamma_ku_k^2)dv_g=\varphi(x_0)+o_k(1).\ee
  \end{lemma}
  {\it Proof}. We only prove the first equality of (\ref{det}), since the proof of the second one is the same.
  Let $(U_{x_0},\psi_{x_0};\{y_1,y_2\})$ be the isothermal coordinate system around $x_0$ given by (\ref{iso-coord}) and (\ref{metr-repre}). We calculate
  by using Lemma \ref{x'} that
  \bna\nonumber
  \int_{\psi_{x_0}^{-1}(\mathbb{B}_{Rr_k}(\widetilde{x}_k)\cap{\mathbb{R}^{2+}})}\varphi\f{c_k }{\lambda_k}|u_k|\exp(\gamma_ku_k^2)dv_g&=&
  \int_{\mathbb{B}_{Rr_k}(\widetilde{x}_k)\cap{\mathbb{R}^{2+}}}{\varphi\circ\psi_{x_0}^{-1}}\f{c_k }{\lambda_k}u_k^\ast\exp(\gamma_k{u_k^\ast}^2)
  \exp(2\psi)dy\\\nonumber
  &=&(\varphi(x_0)+o_k(1))\int_{\mathbb{B}_R\cap\{z_2>-y_{2,k}/r_k\}}\exp(4\pi\eta_0)dz\\
  &=&\varphi(x_0)\int_{\mathbb{B}_R\cap\mathbb{R}^{2+}}\exp(4\pi\eta_0)dz+o_k(1)\\
  &=&\varphi(x_0)+o_k(1)+o_R(1).
  \ena
  Also we have for any fixed $0<\beta<1$,
  \bna
  \int_{\{u_k>\beta c_k\}\setminus\psi_{x_0}^{-1}(\mathbb{B}_{Rr_k}(\widetilde{x}_k)}|\varphi|\f{c_k }{\lambda_k}u_k\exp(\gamma_ku_k^2)dv_g
  &\leq&C\le(1-\int_{\mathbb{B}_{Rr_k}(\widetilde{x}_k)\cap{\mathbb{R}^{2+}}}\f{{u_k^\ast}^2}{\lambda_k}
  \exp(\gamma_k{u_k^\ast}^2)\exp(2f)dy\ri)\\
  &=&C\le(1-\int_{\mathbb{B}_R\cap\mathbb{R}^{2+}}\exp(4\pi\eta_0)dy+o_k(1)\ri)\\
  &=&C(o_k(1)+o_R(1)),
  \ena
  where $C$ is a constant depending only on $\beta$ and $\max_{\overline{\Sigma}}|\varphi|$. Finally we estimate by (\ref{uk-beta})
  and Lemma \ref{concern} that
  \bna
  \int_{u_k\leq \beta c_k}|\varphi|\f{c_k }{\lambda_k}|u_k|\exp(\gamma_ku_k^2)dv_g\leq C\f{c_k}{\lambda_k}\int_{\Sigma}
  |u_{k,\beta}|\exp(\gamma_k u_{k,\beta}^2)dv_g=o_k(1),
  \ena
  where we used $c_k/\lambda_k=o_k(1)$, which is a consequence of Lemma \ref{jixian}.
  Combining the above three estimates, we conclude (\ref{det}). $\hfill\Box$ \\

  The convergence of $c_ku_k$ away from $x_0$ can be described as
  \begin{lemma}\label{Green}
  For any $1<q<2$, $c_ku_k$ converges to $G_{\alpha,x_0}$  weakly in $W^{1,q}(\Sigma)$, strongly in $L^s(\Sigma)$ with $s<2q/(2-q)$, and
  in $C^1_{\rm loc}(\overline{\Sigma}\setminus\{x_0\})$ as $k\ra\infty$, where $G_{\alpha,x_0}$ satisfies
  \be\label{G-x0-1}\le\{\begin{array}{lll}
  \Delta_{g}{G}_{\alpha,x_0}-\alpha G_{\alpha,x_0}=\delta_{x_0}-\f{1}{{\rm Area}(\Sigma)}&{\rm in}&\Sigma\\[1.5ex]
  \f{\p}{\p \bm\nu}{G}_{\alpha,x_0}=0&{\rm on}&\p\Sigma\\[1.5ex]
  \int_\Sigma {G}_{\alpha,x_0}dv_{g}=0
  \end{array}\ri.\ee
  in the distributional sense.
  \end{lemma}
  {\it Proof}. In view of (\ref{E-L}), $c_ku_k$ is a solution of
  \be\label{ckuk}\Delta_g(c_ku_k)-\alpha c_ku_k=f_k\equiv\f{1}{\lambda_k}c_ku_k\exp(\gamma_k u_k^2)-\f{c_k\mu_k}{\lambda_k}\quad{\rm in}\quad \Sigma.\ee
  Integrating both sides of (\ref{ckuk}) and recalling Lemma \ref{tends-delta}, we conclude that
  \be\label{vol}c_k\mu_k/\lambda_k\ra 1/{\rm Area}(\Sigma)\quad{\rm as}\quad  k\ra\infty,\ee and that
  $f_k$ is bounded in $L^1(\Sigma,g)$. We {\it claim} that $c_ku_k$ is also bounded in $L^1(\Sigma,g)$. Suppose not. Let $v_k=c_ku_k/\|c_ku_k\|_1$,
  where $\|\cdot\|_1$ denotes the $L^1(\Sigma,g)$ norm.
  Then $\|v_k\|_1=1$ and satisfies
  \be\label{v-equ}\Delta_g v_k=\alpha v_k+f_k/\|c_ku_k\|_1\quad{\rm in}\quad \Sigma.\ee
  By the Green representation formula (Lemma \ref{G-N}),
  $$v_k(x)=\int_\Sigma \textsf{G}(x,y)\Delta_gv_k(y)dv_{g,y}.$$
  Recalling (\ref{grad}), we have for any $1<q<2$ by using the H\"older inequality and the Fubini theorem,
  \bea\nonumber
  \int_\Sigma|\nabla v_k|^qdv_g&\leq& \int_\Sigma \le(\int_\Sigma|\nabla_{g,x}\textsf{G}(x,y)|^q|\Delta_gv_k(y)|dv_{g,y}\ri)
  \le(\int_\Sigma|\Delta_gv_k(y)|dv_{g,y}\ri)^{q-1}dv_{g,x}\\\nonumber &\leq& \|\Delta_gv_k\|_1^q\sup_{y\in\Sigma}\|\nabla_{g}\textsf{G}(\cdot,y)\|_q^q\\
  &\leq&C.\label{delta}
  \eea
  This together with the Poincare inequality implies that $v_k$ is bounded in $W^{1,q}(\Sigma,g)$. Then up to a subsequence, we assume
  $v_k$ converges to $v_0$ weakly in $W^{1,q}(\Sigma,g)$, strongly in $L^s(\Sigma,g)$ with $s<2q/(2-q)$, and almost everywhere in $\Sigma$.
  As a consequence, $\|v_0\|_1=1$ and $v_0$ is a distributional solution of
  \be\label{v0}\Delta_gv_0-\alpha v_0=0\quad{\rm in}\quad \Sigma,\ee
  where we have used (\ref{v-equ}) and $f_k/\|c_ku_k\|_1\ra 0$  in $L^1(\Sigma,g)$ as $k\ra\infty$. Since $\alpha<\lambda_g(\Sigma)$, it follows
  from (\ref{v0}) that $v_0\equiv 0$ in $\Sigma$, which contradicts $\|v_0\|_1=1$. Hence we conclude our claim $\|c_ku_k\|_1\leq C$. Then coming back
  to (\ref{ckuk}), we see that $\Delta_g(c_ku_k)$ is bounded in $L^1(\Sigma,g)$. In the same way as (\ref{delta}), we obtain
  $$\int_\Sigma |\nabla_g(c_ku_k)|^qdv_g\leq C.$$
  Hence $c_ku_k$ is bounded in $W^{1,q}(\Sigma,g)$. There exists some $G_{\alpha,x_0}$ such that
  $c_ku_k$ converges to $G_{\alpha, x_0}$ weakly in $W^{1,q}(\Sigma,g)$, strongly in $L^s(\Sigma,g)$ for any $s<2q/(2-q)$, and almost everywhere in $\Sigma$.
  In view of Lemma \ref{tends-delta}, $G_{\alpha, x_0}$ satisfies (\ref{G-x0-1})
  in the distributional sense.

  It follows from (\ref{L0}), (\ref{vol}), Lemmas \ref{lemma3} and \ref{concern} that for any $\Omega^\prime\subset\subset\overline{\Sigma}\setminus\{x_0\}$,
  there exists some $p>2$ such that $f_k$ is bounded in $L^p(\Omega^\prime)$. Applying elliptic estimates to (\ref{ckuk}), we have $c_ku_k\ra G_{\alpha, x_0}$ in
  $C^1(\overline{\Omega^{\prime\prime}})$ for any $\Omega^{\prime\prime}\subset\Omega^\prime$. This ends the proof of the lemma. $\hfill\Box$\\

  The function $G_{\alpha,x_0}$ can be decomposed near $x_0$ as below.

  \begin{lemma}\label{G-decom}
  In the isothermal coordinate system (\ref{iso-coord}) around $x_0$, the function $G_{\alpha,x_0}$ can be written as
  the form
  \be\label{G-rep-eucl}G_{\alpha,x_0}\circ\psi_{x_0}^{-1}(y)=-\f{1}{\pi}\log |y|+h(y),\quad \forall y\in\overline{\mathbb{B}_{r_0}^+}\setminus\{0\},\ee
  where $h\in C^1(\overline{\mathbb{B}_{r_0}^+})$.
  \end{lemma}

  {\it Proof}. In the isothermal coordinate system (\ref{iso-coord}) near $x_0$, we set
  $${G}^\ast_{\alpha,x_0}(y)=\le\{\begin{array}{lll}
  G_{\alpha,x_0}\circ\psi_{x_0}^{-1}(y_1,y_2)&{\rm if}&y_2\geq 0\\[1.5ex]
  G_{\alpha,x_0}\circ\psi_{x_0}^{-1}(y_1,-y_2) &{\rm if}&y_2< 0.
  \end{array}\ri.$$
  It follows from (\ref{G-x0-1}) and the fact $\p {G}^\ast_{\alpha,x_0}/\p y_2=0$ on $\p\mathbb{B}_{r_0}^+\cap \p\mathbb{R}^{2+}$ that
  $G^\ast_{\alpha,x_0}$ satisfies
  \be\label{G-dis}-\Delta_{\mathbb{R}^2}G^\ast_{\alpha,x_0}-\alpha\exp(2f)G^\ast_{\alpha,x_0}=2\delta_0-\exp(2f)/{\rm Area}(\Sigma)\quad{\rm in}\quad
  \mathbb{B}_{r_0}\ee
  in the distributional sense. Namely, for any $\varphi\in C_0^2(\mathbb{B}_{r_0})$, there holds
  $$-\int_{\mathbb{B}_{r_0}} (\Delta_{\mathbb{R}^{2}}\varphi)G^\ast_{\alpha,x_0} dy-\alpha\int_{\mathbb{B}_{r_0}}\varphi \exp(2f)
  G^\ast_{\alpha,x_0} dy=
  2\varphi(0)-\int_{\mathbb{B}_{r_0}}\varphi \exp(2f)/{\rm Area}(\Sigma)dy.$$
  Noting also
  \be\label{log}\Delta_{\mathbb{R}^2}\log|y|=2\pi\delta_0\quad{\rm in}\quad
  \mathbb{B}_{r_0},\ee
  we obtain by subtracting (\ref{log}) from (\ref{G-dis}) that
  \be\label{g-reg}-\Delta_{\mathbb{R}^2}\le(G^\ast_{\alpha,x_0}+\f{1}{\pi}\log|y|\ri)=\alpha\exp(2f)G^\ast_{\alpha,x_0}-\exp(2f)/{\rm Area}(\Sigma)
  \quad{\rm in}\quad
  \mathbb{B}_{r_0}.\ee
  Then (\ref{G-rep-eucl}) follows immediately from elliptic estimates on (\ref{g-reg}). $\hfill\Box$\\


 Let $\widetilde{x}_{0,k}$ and $r_0$ be given as in (\ref{x-0k}) and (\ref{iso-coord}) respectively. For any real numbers $R>0$ and $0<s<r_0$, we denote
 $$\mathbb{T}_k^+=\mathbb{B}_{s}^+(\widetilde{x}_{0,k})\setminus\mathbb{B}_{Rr_k}^+(\widetilde{x}_{0,k}),\,\,\,
 \Gamma_{s,k}^+=\p\mathbb{B}_{s}^+(\widetilde{x}_{0,k})\setminus\p\mathbb{R}^{2+},\,\,\,
 \Gamma_{R,k}^+=
 \p\mathbb{B}_{Rr_k}^+(\widetilde{x}_{0,k})\setminus\p\mathbb{R}^{2+}$$
 and
 $$m_{s,k}=\sup_{\Gamma_{s,k}^+}u_k\circ\psi_{x_0}^{-1},\quad i_{R,k}=\inf_{\Gamma_{R,k}^+}u_k\circ\psi_{x_0}^{-1}.$$
 In view of Lemmas \ref{x'}, \ref{Green} and \ref{G-decom}, there holds
 \be\label{inf-sup}\le\{\begin{array}{lll}m_{s,k}=\f{-\f{1}{\pi}\log s+h(0)+o_k(1)+o_s(1)}{c_k}\\[1.5ex]
 i_{R,k}=c_k+\f{-\f{1}{2\pi}\log(1+\f{\pi}{2}R^2)+o_k(1)}{c_k},\end{array}\ri.\ee
 where $o_s(1)\ra 0$ as $s\ra 0+$. Define a sequence of function sets
 $$\mathscr{S}_k=\le\{u\in W^{1,2}(\mathbb{T}_k^+): u|_{\Gamma_{s,k}^+}=m_{s,k},\,u|_{\Gamma_{R,k}^+}=i_{R,k}\ri\}.$$
 Since $i_{R,k}> m_{s,k}$ for sufficiently large $k$, the Poincare inequality holds on $\mathscr{S}_k$. By a direct method of variation,
 it then follows that
 \be\label{inf}\inf_{u\in\mathscr{S}_k}\int_{\mathbb{T}_k^+}|\nabla_{\mathbb{R}^2}u|^2dy\ee
 can be attained by the  harmonic function
 $$h_k(y)=\f{m_{s,k}(\log|y-\widetilde{x}_{0,k}|-\log(Rr_k))+i_{R,k}(\log s-\log|y-\widetilde{x}_{0,k}|)}{\log s-\log(Rr_k)}.$$
 As a consequence
 \be\label{h}\int_{\mathbb{T}_k^+}|\nabla_{\mathbb{R}^2}h|^2dy=\f{\pi(m_{s,k}-i_{R,k})^2}{\log s-\log(Rr_k)}.\ee
 Define a sequence of functions
 $$\widetilde{u}_k(y)=\max\{m_{s,k},\min\{u_k\circ \psi_{x_0}^{-1}(y),i_{R,k}\}\},\quad y\in\mathbb{B}_{r_0}^+.$$
 One can see that $\widetilde{u}_k$ belongs to $\mathscr{S}_k$ and that
 \bea\nonumber\int_{\mathbb{T}_k^+}|\nabla_{\mathbb{R}^{2}}\widetilde{u}_k|^2dy&\leq&
 \int_{\mathbb{T}_k^+}|\nabla_{\mathbb{R}^{2}}({u}_k\circ\psi_{x_0}^{-1})|^2dy\\\nonumber
 &=&\int_{\psi_{x_0}^{-1}(\mathbb{T}_k^+)}|\nabla_gu_k|^2dv_g\\\nonumber
 &=&1+\alpha\int_\Sigma u_k^2dv_g-\int_{\Sigma\setminus\psi_{x_0}^{-1}(\mathbb{B}_{s}^+(\widetilde{x}_{0,k}))}|\nabla_gu_k|^2dv_g\\
 &&-\int_{\mathbb{B}_{Rr_k}^+(\widetilde{x}_{0,k})}|\nabla_{\mathbb{R}^2}(u_k\circ\psi_{x_0}^{-1})|^2dy.\label{ener}\eea
 Combining (\ref{inf}), (\ref{h}) and (\ref{ener}), we have
 \bea\f{\pi(m_{s,k}-i_{R,k})^2}{\log s-\log(Rr_k)}&\leq&
 1+\alpha\int_\Sigma u_k^2dv_g-\int_{\Sigma\setminus\psi_{x_0}^{-1}(\mathbb{B}_{s}^+(\widetilde{x}_{0,k}))}|\nabla_gu_k|^2dv_g\nonumber\\
 &&-\int_{\mathbb{B}_{Rr_k}^+(\widetilde{x}_{0,k})}|\nabla_{\mathbb{R}^2}(u_k\circ\psi_{x_0}^{-1})|^2dy.\label{inequ}\eea
 It follows from (\ref{inf-sup}) that
 \be\label{left}\f{\pi(m_{s,k}-i_{R,k})^2}{\log s-\log(Rr_k)}=\f{2\pi c_k^2-2\log(1+\f{\pi}{2}R^2)+4\log s-4\pi h(0)+o(1)}
 {\gamma_k c_k^2+2\log s-2\log R-\log\f{\lambda_k}{c_k^2}}.\ee
 Let $\bm\nu$ be the unit outward vector on $\p\psi_{x_0}^{-1}(\mathbb{B}_{s}^+(\widetilde{x}_{0,k}))$.
  We write $\bm\nu=\nu^1\p/\p y_1+\nu^2\p/\p y_2$. Then there holds on $\p\psi_{x_0}^{-1}(\mathbb{B}_{s}^+(\widetilde{x}_{0,k}))\setminus\p\Sigma$,
  \bna\f{\p G_{\alpha,x_0}}{\p\bm\nu}&=&\nu^1\f{\p}{\p y_1}(G_{\alpha,x_0}\circ \psi_{x_0}^{-1})+\nu^2\f{\p}{\p y_2}(G_{\alpha,x_0}\circ \psi_{x_0}^{-1}),\\[1.2ex]
  1=|\bm\nu|^2&=&
  \exp(2f(y))((\nu^1)^2+(\nu^2)^2),\\[1.2ex] d\sigma_g&=&\exp(f(y))d\sigma_0,\ena
  and
  \bna\nonumber
  \int_{\p\psi_{x_0}^{-1}(\mathbb{B}_{s}^+(\widetilde{x}_{0,k}))\setminus\p\Sigma}G_{\alpha,x_0}\f{\p G_{\alpha,x_0}}{\p\bm\nu}d\sigma_g&=&
  \int_{\p\mathbb{B}_{s}^+(\widetilde{x}_{0,k})\setminus\p\mathbb{R}^{2+}}(G_{\alpha,x_0}\circ\psi_{x_0}^{-1})
  \f{\p}{\p\bm\nu_0}(G_{\alpha,x_0}\circ \psi_{x_0}^{-1})
  \exp(f)d\sigma_0\\[1.2ex]\nonumber&=&\int_0^\pi \le(-\f{\log s}{\pi}+h(s\cos t,s\,\sin t)\ri)\le(-\f{1}{\pi s}+\f{\p h}{\p s}
  \ri)sdt+o_k(1)\\[1.2ex]\label{gree}
  &=&-\f{1}{\pi}\log s+h(0)+o_s(1)+o_k(1),
  \ena
  where $\bm\nu_0=(\nu^1,\nu^2)$ is a normal vector field on $\p\mathbb{B}_{s}^+(\widetilde{x}_{0,k})\setminus\p\mathbb{R}^{2+}$
  and $d\sigma_0$ denotes its Euclidean arc length element. Thus
  \bna
  \int_{\Sigma\setminus\psi_{x_0}^{-1}(\mathbb{B}_{s}^+(\widetilde{x}_{0,k}))}|\nabla_gG_{\alpha,x_0}|^2dv_g&=&
  \int_{\Sigma\setminus\psi_{x_0}^{-1}(\mathbb{B}_{s}^+(\widetilde{x}_{0,k}))}G_{\alpha,x_0}\Delta_gG_{\alpha,x_0} dv_g\\[1.2ex]&&+
  \int_{\p\psi_{x_0}^{-1}(\mathbb{B}_{s}^+(\widetilde{x}_{0,k}))\setminus\p\Sigma}G_{\alpha,x_0}\f{\p G_{\alpha,x_0}}{\p\bm\nu}d\sigma_g\\[1.2ex]
  &=&\alpha\int_\Sigma G_{\alpha,x_0}^2dv_g-\f{1}{\pi}\log s+A_0+o_k(1)+o_s(1).
  \ena
  This together with Lemma \ref{Green} leads to
  \bea\nonumber
  \int_{\Sigma\setminus\psi_{x_0}^{-1}(\mathbb{B}_{s}^+(\widetilde{x}_{0,k}))}|\nabla_gu_k|^2dv_g&=&\f{1}{c_k^2}\le(
  \int_{\Sigma\setminus\psi_{x_0}^{-1}(\mathbb{B}_{s}^+(\widetilde{x}_{0,k}))}|\nabla_gG_{\alpha,x_0}|^2dv_g+o_k(1)\ri)\\
  &=&\f{1}{c_k^2}\le(\alpha\int_\Sigma G_{\alpha,x_0}^2dv_g-\f{1}{\pi}\log s+h(0)+o_k(1)+o_s(1)\ri).\label{33}
  \eea
  By Lemma \ref{x'}, we obtain
  \bea\nonumber
  \int_{\psi_{x_0}^{-1}(\mathbb{B}_{Rr_k}^+(\widetilde{x}_{0,k}))}|\nabla_gu_k|^2dv_g&=&
  \int_{\mathbb{B}_{Rr_k}^+(\widetilde{x}_{0,k}))}|\nabla_{\mathbb{R}^2}(u_k\circ \psi_{x_0}^{-1})|^2dy\\&=&
  \f{1}{c_k^2}
  \le(\f{1}{2\pi}\log(1+\f{\pi}{2}R^2)-\f{1}{2\pi}+o_k(1)+o_R(1)\ri),  \label{inner-energy}
  \eea
  where $o_R(1)\ra 0$ as $R\ra\infty$.
  Combining (\ref{inequ}), (\ref{left}), (\ref{33}), (\ref{inner-energy}) and passing to the limit
  $k\ra\infty$ firstly, then $R\ra\infty$ and finally $s\ra 0$, we calculate
  $$\limsup_{k\ra\infty}\f{\lambda_k}{c_k^2}\leq \f{\pi}{2}\exp(1+2\pi h(0)),$$
  which together with (\ref{L1}) and (\ref{sup-2pi}) leads to
  \be\label{76'}\sup_{u\in\mathcal{H},\|u\|_{1,\alpha}\leq 1}\int_\Sigma \exp(2\pi u^2)dv_g=\lim_{k\ra\infty}\int_\Sigma \exp(\gamma_k u_k^2)dv_g\leq {\rm Area}(\Sigma)+\f{\pi}{2}\exp(1+2\pi h(0)).\ee

  \subsection{Test function computation}\label{4.4}
  We shall construct a sequence of functions $\phi_k\in\mathcal{H}$ with $\|\phi_k\|_{1,\alpha}=1$ such that
  \be\label{gret}\int_\Sigma \exp(2\pi \phi_k^2)dv_g>{\rm Area}(\Sigma)+\f{\pi}{2}\exp(1+2\pi h(0)).\ee
  The contradiction between (\ref{gret}) and (\ref{76'}) implies that (\ref{c-}) can not hold. Then applying elliptic estimates to (\ref{E-L}),
  we finish the proof of Theorem \ref{thm1}.

  To proceed, we use the isothermal coordinate system $(U_{x_0},\psi_{x_0};\{y_1,y_2\})$, which is defined as in (\ref{iso-coord}), and let
  $$\widetilde{\phi}_k(y)=\le\{\begin{array}{lll}c+\f{1}{c}\le(-\f{1}{2\pi}\log(1+\f{\pi}{2}|ky|^2)+B\ri)&{\rm when}& |y|\leq \f{\log k}{k}
  \\[1.2ex] \f{1}{c}(G_{\alpha,x_0}\circ\psi_{x_0}^{-1}(y)-\eta(y)\beta(y)) &{\rm when}&\f{\log k}{k}<|y|<2\f{\log k}{k}\end{array}\ri.$$
  where $\beta(y)=G_{\alpha,x_0}\circ\psi_{x_0}^{-1}(y)+\f{1}{\pi}\log|y|-h(0)$, $\eta(y)=\eta(|y|)$ is a radially symmetric function satisfying
  $\eta\in C_0^1(\mathbb{B}_{2k^{-1}{\log k}})$,
  $\eta\equiv 1$ in $\mathbb{B}_{k^{-1}{\log k}}$, $\|\nabla_{\mathbb{R}^2}\eta\|_{L^\infty}=O(\f{k}{\log k})$, $B$ and $c$ are constants
  depending only on $k$
  to be determined later. Define
  \be\label{phi-k}\phi_k=\le\{\begin{array}{lll}
  \widetilde{\phi}_k\circ \psi_{x_0}&{\rm on}& \psi_{x_0}^{-1}(\mathbb{B}_{2k^{-1}{\log k}}^+)\\[1.2ex]
  \f{G_{\alpha,x_0}}{c}&{\rm on}& \Sigma\setminus\psi_{x_0}^{-1}(\mathbb{B}_{2k^{-1}{\log k}}^+).
  \end{array}\ri.\ee
  On $\p\mathbb{B}_{k^{-1}{\log k}}^+\setminus\mathbb{R}^{2+}$, we let
  \be\label{c+c}c+\f{1}{c}\le(-\f{1}{2\pi}\log(1+\f{\pi}{2}(\log k)^2)+B\ri)=\f{1}{c}\le(-\f{1}{\pi}\log(k^{-1}\log k)+h(0)\ri),\ee
  which leads to $\phi_k\in W^{1,2}(\Sigma,g)$. It follows from (\ref{c+c}) that
  \be\label{c2}2\pi c^2=2\log k-2\pi B+2\pi h(0)+\log\f{\pi}{2}+O\le(\f{1}{(\log k)^2}\ri).\ee
  Clearly we calculate
  \bea\nonumber
  \int_{\psi_{x_0}^{-1}(\mathbb{B}_{k^{-1}\log k}^+)}|\nabla_g\phi_k|^2dv_g&=&\int_{\mathbb{B}_{k^{-1}\log k}^+}
  |\nabla_{\mathbb{R}^2}\widetilde{\phi}_k|^2dy\\ \nonumber
  &=&\f{1}{4 c^2}\int_{\mathbb{B}_{\log k}^+}\f{|z|^2}{(1+\f{\pi}{2}|z|^2)^2}dz\\\label{energy-1}
  &=&\f{1}{2\pi c^2}\le(2\log(\log k)+\log\f{\pi}{2}-1+O\le(\f{1}{(\log k)^2}\ri)\ri).
  \eea
  Denoting $\mathbb{T}_k^+=\mathbb{B}_{2k^{-1}\log k}^+\setminus\mathbb{B}_{k^{-1}\log k}^+$, we have
  \bea\nonumber
  \int_{\psi_{x_0}^{-1}(\mathbb{T}_k^+)}|\nabla_g\phi_k|^2dv_g&=&\int_{\mathbb{T}_k^+}|\nabla_{\mathbb{R}^{2}}\widetilde{\phi}_k|^2dy\\
  \nonumber&=&\int_{\mathbb{T}_k^+}\f{1}{c^2}|\nabla_{\mathbb{R}^2}(G_{\alpha,x_0}\circ\psi_{x_0}^{-1})|^2dy+\int_{\mathbb{T}_k^+}\f{1}{c^2}
  |\nabla_{\mathbb{R}^2}(\eta \beta)|^2dy\\\nonumber
  &&\quad-\int_{\mathbb{T}_k^+}\f{2}{c^2}\nabla_{\mathbb{R}^2}(G_{\alpha,x_0}\circ\psi_{x_0}^{-1})\nabla_{\mathbb{R}^2}(\eta \beta)dy\\
  &=&\f{1}{c^2}\le(\f{2}{\pi}\log 2+O\le(\f{1}{(\log k)^2}\ri)\ri).\label{energy-2}
  \eea
    Writing $\widetilde{G}= G_{\alpha,x_0}\circ\psi_{x_0}^{-1}$ and $\bm\nu=v^1\p/\p y_1+v^2\p/\p y_2$, we get on
    $\psi_{x_0}^{-1}(\mathbb{B}_{2k^{-1}\log k}^+)\cap\p\Sigma$,
  $$\f{\p G_{\alpha,x_0}}{\p\bm{\nu}}=v^1\f{\p\widetilde{G}}{\p y_1}+v^2\f{\p\widetilde{G}}{\p y_2}=\exp(-f)\f{\p\widetilde{G}}{\p \bm\nu_0},$$
  where $\bm\nu_0=(\psi_{x_0})_\ast(\bm\nu)/|(\psi_{x_0})_\ast(\bm\nu)|$ is the unit outward vector field on $\p\mathbb{B}_{2k^{-1}\log k}^+
  \setminus\p\mathbb{R}^{2+}$. Moreover
  $d\sigma_g=\exp(f)d\sigma_0$, where $d\sigma_0$ is the Euclidean arc-length element of $\p\mathbb{B}_{2k^{-1}\log k}^+
  \setminus\p\mathbb{R}^{2+}$. It then follows that
  \bna\nonumber
  \int_{\psi_{x_0}^{-1}(\p\mathbb{B}_{2k^{-1}\log k}^+
  \setminus\p\mathbb{R}^{2+})}G_{\alpha,x_0}\f{\p G_{\alpha,x_0}}{\p\bm\nu}d\sigma_g&=&\int_{\p\mathbb{B}_{2k^{-1}\log k}^+
  \setminus\p\mathbb{R}^{2+}}\widetilde{G}\f{\p \widetilde{G}}{\p\bm\nu_0}d\sigma_0\\\nonumber
  &=&\int_{\p\mathbb{B}_{2k^{-1}\log k}^+\setminus\p\mathbb{R}^{2+}}\le(-\f{1}{\pi}\log|y|+h(0)+O(|y|)\ri)\\\nonumber
  &&\quad\times\le(-\f{1}{\pi|y|}+O(1)\ri)d\sigma_0\\
  &=&\f{\log 2}{\pi}+\f{1}{\pi}\log\le(\f{\log k}{k}\ri)-h(0)+O\le(\f{1}{(\log k)^2}\ri).\label{arc}
  \ena
  This together with
  $$\int_{\Sigma\setminus\psi_{x_0}^{-1}(\mathbb{B}_{2k^{-1}\log k}^+)}G_{\alpha,x_0}^2dv_g=\int_\Sigma G_{\alpha,x_0}^2dv_g+
  O\le(\f{1}{(\log k)^2}\ri)$$
  and
  $$\int_{\Sigma\setminus\psi_{x_0}^{-1}(\mathbb{B}_{2k^{-1}\log k}^+)}G_{\alpha,x_0}dv_g=-\int_{\psi_{x_0}^{-1}(\mathbb{B}_{2k^{-1}\log k}^+)}
  G_{\alpha,x_0}dv_g=O\le(\f{1}{(\log k)^2}\ri)$$
  leads to
  \bea\nonumber
  \int_{\Sigma\setminus\psi_{x_0}^{-1}(\mathbb{B}_{2k^{-1}\log k}^+)}|\nabla_g\phi_k|^2dv_g&=&
  \int_{\Sigma\setminus\psi_{x_0}^{-1}(\mathbb{B}_{2k^{-1}\log k}^+)}\f{|\nabla_gG_{\alpha,x_0}|^2}{c^2}dv_g\\\nonumber
  &=&\f{1}{c^2}\int_{\p(\Sigma\setminus\psi_{x_0}^{-1}(\mathbb{B}_{2k^{-1}\log k}^+))}G_{\alpha,x_0}
  \f{\p G_{\alpha,x_0}}{\p\bm\nu}d\sigma_g\\\nonumber
  &&+\f{1}{c^2}\int_{\Sigma\setminus\psi_{x_0}^{-1}(\mathbb{B}_{2k^{-1}\log k}^+)}G_{\alpha,x_0}\Delta_gG_{\alpha,x_0}dv_g\\\nonumber
  &=&\f{1}{c^2}\le\{\int_{\psi_{x_0}^{-1}(\p\mathbb{B}_{2k^{-1}\log k}^+
  \setminus\p\mathbb{R}^{2+})}G_{\alpha,x_0}\f{\p G_{\alpha,x_0}}{\p\bm\nu}d\sigma_g
  +\alpha\int_{\Sigma}G_{\alpha,x_0}^2dv_g\ri.\\\nonumber
  &&\quad\le.+\f{1}{{\rm Area}(\Sigma)}\int_{
  \psi_{x_0}^{-1}(\mathbb{B}_{2k^{-1}\log k}^+)}G_{\alpha,x_0}dv_g\ri\}\\\nonumber
  &=&\f{1}{c^2}\le\{-\f{1}{\pi}\log\le(\f{\log k}{k}\ri)-\f{\log 2}{\pi}+h(0)+\alpha\int_\Sigma G_{\alpha,x_0}^2dv_g\ri.\\
  &&\qquad\le.
  +O\le(\f{1}{(\log k)^2}\ri)\ri\}.\label{energy-3}
  \eea
  Combining (\ref{energy-1}), (\ref{energy-2}) and (\ref{energy-3}), we conclude
  \bea\label{en-sigma}
  \int_\Sigma|\nabla_g\phi_k|^2dv_g=\f{1}{c^2}\le\{\f{\log k}{\pi}+h(0)+\f{1}{2\pi}\log\f{\pi}{2}-
  \f{1}{2\pi}+\alpha\int_\Sigma G_{\alpha,x_0}^2dv_g+O\le(\f{1}{(\log k)^2}\ri)\ri\}.
  \eea
  Also one can compute
  $$\overline{\phi}_k=\f{1}{{\rm Area}(\Sigma)}\int_\Sigma\phi_kdv_g=\f{1}{c}O\le(\f{1}{(\log k)^2}\ri)$$
  and
  $$\int_\Sigma(\phi_k-\overline{\phi}_k)^2dv_g=\f{1}{c^2}\le(\int_\Sigma G_{\alpha,x_0}^2dv_g+O\le(\f{1}{(\log k)^2}\ri)\ri).$$
  This together with (\ref{en-sigma}) gives
  \bea\nonumber
  \|\phi_k-\overline{\phi}_k\|_{1,\alpha}^2&=&\int_\Sigma|\nabla_g(\phi_k-\overline{\phi}_k)|^2dv_g-\alpha\int_\Sigma(\phi_k-\overline{\phi}_k)^2dv_g\\
  &=&\f{1}{c^2}\le\{\f{\log k}{\pi}+h(0)+\f{1}{2\pi}\log\f{\pi}{2}-
  \f{1}{2\pi}+O\le(\f{1}{(\log k)^2}\ri)\ri\}.\label{norm=1}
  \eea
  Now we set
  \be\label{111}\|\phi_k-\overline{\phi}_k\|_{1,\alpha}=1.\ee
  It follows from (\ref{norm=1}) and (\ref{111}) that
  \be\label{c-22}c^2=\f{\log k}{\pi}+h(0)+\f{1}{2\pi}\log\f{\pi}{2}-
  \f{1}{2\pi}+O\le(\f{1}{(\log k)^2}\ri).\ee
  Inserting (\ref{c-22}) into (\ref{c2}), we obtain
  \be\label{B}B=\f{1}{2\pi}+O\le(\f{1}{(\log k)^2}\ri).\ee
  In view of (\ref{phi-k}), (\ref{c-22}) and (\ref{B}), there holds
  \bna
  \int_{\psi_{x_0}^{-1}(\mathbb{B}_{k^{-1}\log k}^+)}\exp(2\pi(\phi_k-\overline{\phi}_k)^2)dv_g&=&
  \int_{\mathbb{B}_{k^{-1}\log k}^+}\exp\le(2\pi(\widetilde{\phi}_k(y)-\overline{\phi}_k)^2+2f(y)\ri)dy\\
  &=&(1+O((\log k)^{-2}))\int_{\mathbb{B}_{k^{-1}\log k}^+}\exp(2\pi\widetilde{\phi}_k^2(y))dy\\
  &\geq&(1+O((\log k)^{-2}))\int_{\mathbb{B}_{\log k}^+}\exp\le\{2\pi c^2-2\log(1+\f{\pi}{2}|z|^2)\ri.\\
  &&\qquad\qquad+4\pi B\le\}\f{1}{k^2}dz\ri.\\
  &=&(1+O((\log k)^{-2}))\f{\pi}{2}\exp(1+2\pi h(0)),
  \ena
  and
  \bna
  \int_{\Sigma\setminus\psi_{x_0}^{-1}(\mathbb{B}_{k^{-1}\log k}^+)}\exp(2\pi(\phi_k-\overline{\phi}_k)^2)dv_g
  &\geq&\int_{\Sigma\setminus\psi_{x_0}^{-1}(\mathbb{B}_{k^{-1}\log k}^+)}(1+2\pi(\phi_k-\overline{\phi}_k)^2)dv_g\\
  &=&{\rm Area}(\Sigma)+\f{2\pi}{c^2}\int_\Sigma {G_{\alpha,x_0}^2}dv_g+O((\log k)^{-2}).
  \ena
  Therefore
  \be\label{ggg}\int_{\Sigma}\exp(2\pi(\phi_k-\overline{\phi}_k)^2)dv_g\geq {\rm Area}(\Sigma)+\f{\pi}{2}\exp(1+2\pi h(0))+\f{2\pi}{c^2}\int_\Sigma {G_{\alpha,x_0}^2}dv_g+O((\log k)^{-2}).\ee
  Since $(\log k)^{-2}=o(c^{-2})$, we have by (\ref{ggg}) that
  $$\int_{\Sigma}\exp(2\pi(\phi_k-\overline{\phi}_k)^2)dv_g>{\rm Area}(\Sigma)+\f{\pi}{2}\exp(1+2\pi h(0))$$
  for sufficiently large $k$. Therefore $\phi_k-\overline{\phi}_k\in\mathcal{H}$ satisfies (\ref{gret}) provided that $k$
    is chosen sufficiently large, and thus the proof of Theorem \ref{thm1}
  is completely finished.

  \section{Proof of Theorem \ref{thm2}}\label{Sec-5}

  In this section, we shall prove Theorem \ref{thm2} by using the same method of proving Theorem \ref{thm1}.
  We only give its outline but emphasize their differences.

  \subsection{The best constant}
  Let $\tau>0$ be a fixed positive real number, $u$ be any function in $W^{1,2}(\Sigma,g)$, $\|u\|_{1,\tau}$ be defined as in (\ref{ineq-2})
  and $\overline{u}=\f{1}{{\rm Area}(\Sigma)}\int_\Sigma udv_g$.
  By the H\"older inequality,
  \be\label{mean-u}\overline{u}^2\leq \f{1}{{\rm Area}(\Sigma)}\int_\Sigma u^2dv_g\leq \f{\|u\|_{1,\tau}}{\tau{\rm Area}(\Sigma)}.\ee
  Hence, if $\|u\|_{1,\tau}\leq 1$, then $\int_\Sigma|\nabla_g(u-\overline{u})|^2dv_g\leq 1$, and the Young inequality together with
  (\ref{mean-u}) implies that
  for any $\epsilon>0$, there holds
  a constant $C$ depending only on $(\Sigma,g)$, $\alpha$ and $\epsilon$ such that
  \be\label{in-3}\int_\Sigma \exp(\alpha u^2)dv_g\leq C\le(\int_\Sigma \exp(\alpha(1+\epsilon)(u-\overline{u})^2)dv_g\ri)^{{1}/{(1+\epsilon)}}.\ee
  Define
  $$\alpha^\ast=\sup\le\{\alpha: \sup_{u\in W^{1,2}(\Sigma,g),\,\|u\|_{1,\tau}\leq 1}\int_\Sigma \exp(\alpha u^2)dv_g<\infty\ri\}.$$  It follows from (\ref{in-3}) and Lemma \ref{lemma3} that
  \be\label{alp-geq}\alpha^\ast\geq 2\pi.\ee

  Let $M_k$ be defined as in (\ref{Moser-sequence}). Then we have
  $$\|M_k\|_{1,\tau,r}^2=\int_{\mathbb{B}_r}(|\nabla_{\mathbb{R}^2}M_k|^2+\tau M_k^2)dy=1+o_k(1)+o_r(1),$$
  where $o_r(1)\ra 0$ as $r\ra 0$.
  For any $\gamma>2\pi$, there holds
  \bna
  \int_{\mathbb{B}_r^+}\exp(\gamma M_k^2/\|M_k\|_{1,\tau,r}^2)dy&\geq&\int_{\mathbb{B}_{rk^{-1/4}}^+}
  \exp(\gamma(1+o_k(1)+o_r(1))M_k^2)dy\\
  &=&\exp\le(\gamma(1+o_k(1)+o_r(1))\f{\log k}{4\pi}\ri)\f{\pi}{2}r^2k^{-1/2}\\
  &=&\f{\pi}{2}r^2k^{\f{\gamma}{4\pi}(1+o_k(1)+o_r(1))-\f{1}{2}}.
  \ena
  Let $(U_{x_0},\psi_{x_0};\{y_1,y_2\})$ be the isothermal coordinate system around $x_0\in\p\Sigma$, and
   the metric $g$ can be written as $g=\exp(2f(y))(dy_1^2+dy_2^2)$. Define a sequence of functoions $\widetilde{M}_k=M_{k,r}\circ \psi_{x_0}$. Then  we have
  $$\|\widetilde{M}_k\|_{1,\tau}^2=\int_\Sigma(|\nabla_g\widetilde{M}_k|^2+\tau \widetilde{M}_k^2)dv_g=1+o_k(1)+o_r(1).$$
  It follows that for any fixed $\gamma>2\pi$, if $r>0$ is chosen sufficiently small,
  \bna
  \int_\Sigma\exp(\gamma \widetilde{M}_k^2/\|\widetilde{M}_k\|_{1,\tau}^2)dv_g&\geq&\int_{U_p}\exp(\gamma \widetilde{M}_k^2/\|\widetilde{M}_k\|_{1,\tau}^2)dv_g\\
  &\geq&\int_{\mathbb{B}_{rk^{-1/4}}^+}\exp(\gamma(1+o_k(1)+o_r(1))M_k^2)\exp(2f)dy\\
  &=&(1+o(1))\pi r^2k^{\f{\gamma}{4\pi}(1+o(1))-\f{1}{2}}\\
  &\ra& +\infty
  \ena
  as $k\ra\infty$. This leads to $\alpha^\ast\leq 2\pi$, which together with (\ref{alp-geq}) implies that
  $\alpha^\ast=2\pi$.

  \subsection{The existence of extremals for the supremums in (\ref{ineq-2})}
  By a direct method of variation, for any $k\in\mathbb{N}$, there exists a nonnegative function $u_k$ with $\|u_k\|_{1,\tau}=1$ such that
  $$\label{maximizer-2}\int_\Sigma \exp(\gamma_ku^2_k)dv_g=\sup_{\|u\|_{1,\tau}\leq 1}\int_\Sigma\exp(\gamma_ku^2)dv_g,$$
  where $\gamma_k=2\pi-1/k$. One can easily check that $u_k$ satisfies the Euler-Lagrange equation
  \be\label{E-L-2}\le\{\begin{array}{lll}
  \Delta_gu_k+\tau u_k=\f{1}{\lambda_k}u_k\exp(\gamma_ku_k^2)&{\rm in}& \Sigma\\[1.2ex]
  u_k>0&{\rm in}&\Sigma\\[1.2ex]
  \f{\p u_k}{\p\bm\nu}=0&{\rm on}&\p\Sigma\\[1.2ex]
  \lambda_k=\int_\Sigma u_k^2\exp(\gamma_ku_k^2)dv_g.
  \end{array}\ri.\ee
  With no loss of generality, we assume $c_k=u_k(x_k)=\max_{\overline{\Sigma}}u_k\ra+\infty$ and $x_k\ra x_0\in\overline{\Sigma}$
  as $k\ra\infty$. Then as in Lemma \ref{concern}, we have $x_0\in\p\Sigma$, $u_k$ converges to $0$ weakly in $W^{1,2}(\Sigma,g)$,
  strongly in $L^q(\Sigma)$ for any $q>1$, and $|\nabla_gu_k|dv_g\rightharpoonup \delta_{x_0}$ in the sense of measure.

   In an isothermal coordinate system $(U_{x_0},\psi_{x_0};\{y_1,y_2\})$ around $x_0$,
  $\psi_{x_0} (U_{x_0})=\mathbb{B}_{r_0}$,
  the metric $g$ can be written as
  $g=\exp(2f(y))(dy_1^2+dy_2^2)$ with $f\in C^1(\overline{\mathbb{B}_{r_0}})$ and $f(0,0)=0$; moreover,
   the unit outward vector field $\bm\nu$ on the boundary
  $\p\Sigma$ can be written as
  $\bm\nu=\exp(-f(y)){\p}/{\p y_2}$. For any $u\in C^1(\overline{\Sigma})$, the normal derivative $\p u/\p\bm{\nu}$ can be represented by
  $$\label{norm-2}\f{\p u}{\p\bm\nu}=\exp(-f(y))\f{\p}{\p y_2}(u\circ\psi_{x_0}^{-1}).$$
  Denote $\widetilde{x}_k=\psi_{x_0}(x_k)=(y_{1,k},y_{2,k})$ and $\widetilde{x}_{0,k}=(y_{1,k},0)$. Let $r_k>0$ satisfy
  $$r_k^2=\f{\lambda_k}{c_k^2}\exp(-\gamma_kc_k^2).$$
  Using the same argument in the proof of Lemmas \ref{bubble} and \ref{x'}, we have as $k\ra\infty$,
  $$\label{converge}c_k(u_k\circ\psi_{x_0}^{-1}(\widetilde{x}_{0,k}+r_k\cdot)-c_k)\ra -\f{1}{2\pi}\log(1+\f{\pi}{2}|\cdot|^2)\quad
  {\rm in}\quad C^1_{\rm loc}(\mathbb{R}^{2+}\cup\p\mathbb{R}^{2+}).$$
   Similar to Lemma \ref{tends-delta}, we also have that
  for any $\varphi\in C^2(\overline{\Sigma})$, there holds
  $$\label{det-2}\int_\Sigma\varphi\f{c_k }{\lambda_k}u_k\exp(\gamma_ku_k^2)dv_g=\varphi(x_0)+o_k(1).$$
  In particular,
  \be\label{L-bdd}\f{1}{\lambda_k}c_k\|u_k\exp(\gamma_ku_k^2)\|_{L^1(\Sigma,g)}\leq C\ee
  and in the sense of measure
  $$\label{ra-delta}\f{1}{\lambda_k}c_ku_k\exp(\gamma_ku_k^2)\rightharpoonup \delta_{x_0}.$$
  In view of (\ref{E-L-2}), there holds
  \be\label{c-u}\Delta_g(c_ku_k)+\tau(c_ku_k)=\f{1}{\lambda_k}c_ku_k\exp(\gamma_k u_k^2)\quad{\rm in}\quad \overline{\Sigma}.\ee
  Integrating both sides of (\ref{c-u}), we have by noticing (\ref{L-bdd}), $u_k> 0$ in $\Sigma$ and $\p u_k/\p\bm\nu=0$ on $\p\Sigma$ that
  \be\label{l-1}\int_\Sigma c_ku_kdv_g\leq C\ee
  and
  $$\label{l-2}\int_\Sigma|\Delta_g(c_ku_k)|dv_g\leq C.$$
  Let
  $$w_k=c_ku_k-\f{1}{{\rm Area}(\Sigma)}\int_\Sigma c_ku_kdv_g.$$
  Then we obtain by using the Green representation formula,
  $$w_k(x)=\int_\Sigma \textsf{G}(x,y)\Delta_gw_k(y)dv_{g,y},$$
  where $\textsf{G}(\cdot,\cdot)$ is defined as in Lemma \ref{G-N}.
  An obvious analog of (\ref{delta}) reads
  $\|\nabla_gw_k\|_{L^q(\Sigma,g)}\leq C$ for all $1<q<2$.
  Hence
  $\|\nabla_g(c_ku_k)\|_{L^q(\Sigma,g)}\leq C$  for all $1<q<2$.
  This together with (\ref{l-1}) implies that $c_ku_k$ is bounded in $W^{1,q}(\Sigma,g)$ for any $1<q<2$.
  Similar to Lemma \ref{Green}, $c_ku_k$ converges to $G_{\tau,x_0}$  weakly in $W^{1,q}(\Sigma,g)$, strongly in $L^s(\Sigma,g)$ with $s<2q/(2-q)$, and
  in $C^1_{\rm loc}(\overline{\Sigma}\setminus\{x_0\})$ as $k\ra\infty$, where $G_{\tau,x_0}$ satisfies in the distributional sense
  $$\label{G-x0}\le\{\begin{array}{lll}
  \Delta_{g}{G}_{\tau,x_0}+\tau G_{\tau,x_0}=\delta_{x_0}&{\rm in}&\Sigma\\[1.5ex]
  \f{\p}{\p \bm\nu}{G}_{\tau,x_0}=0&{\rm on}&\p\Sigma\\[1.5ex]
  \int_\Sigma {G}_{\tau,x_0}dv_{g}=0.
  \end{array}\ri.$$
    Similar to Lemma \ref{G-decom}, in the isothermal coordinate system
  $(U_{x_0},\psi_{x_0};\{y_1,y_2\})$ near $x_0$, we have
  $$G_{\tau,x_0}\circ\psi_{x_0}^{-1}(y)=-\f{1}{\pi}\log|y|+h(y),$$
  where $h\in C^1(\overline{\mathbb{B}_{r_0}^+})$. Then repeating the argument
  of deriving (\ref{76'}), we obtain
  \be\label{gge-0}
  \sup_{\|u\|_{1,\tau}\leq 1}\int_\Sigma\exp(2\pi u^2)dv_g=\lim_{k\ra\infty}\int_\Sigma\exp(\gamma_ku_k^2)dv_g
  \leq{\rm Area}(\Sigma)+\f{\pi}{2}\exp(1+2\pi h(0)).
  \ee

  Let $\phi_k$ be defined as in (\ref{phi-k}). We first require $\phi_k\in W^{1,2}(\Sigma,g)$. In particular, (\ref{c+c}) and
  thus (\ref{c2}) hold. A straightforward calculation shows
  $$\int_\Sigma(|\nabla_g\phi_k|^2+\tau \phi_k^2)dv_g=\f{1}{c^2}\le\{\f{\log k}{\pi}+h(0)+\f{1}{2\pi}\log\f{\pi}{2}-
  \f{1}{2\pi}+O\le(\f{1}{(\log k)^2}\ri)\ri\}.$$
  We further require
  $$\|\phi_k\|_{1,\tau}^2=\int_\Sigma(|\nabla_g\phi_k|^2+\tau \phi_k^2)dv_g=1.$$
  It then follows that (\ref{c-22}) and (\ref{B}) hold. As a consequence, we calculate as before
  \be\label{gge}\int_\Sigma\exp(2\pi\phi_k^2)dv_g>{\rm Area}(\Sigma)+\f{\pi}{2}\exp(1+2\pi h(0)),\ee
  provided that $k$ is sufficiently large.

  The contradiction between (\ref{gge-0}) and (\ref{gge}) implies that $c_k$ must be bounded, and thus the supremum in (\ref{ineq-2})
  can be attained for $\gamma=2\pi$. This ends the proof of Theorem \ref{thm2}.

  \section{Proof of Theorem \ref{thm3}}\label{Sec-6}

  {\it Proof of Theorem \ref{thm3}}. Let $0\leq\alpha<\lambda_{\rm N}(\Sigma)$ and $\tau>0$ be two fixed real numbers.

  {\it The inequality (\ref{in-01}) implies the inequality (\ref{in-02})}. Suppose (\ref{in-01}) holds. To derive (\ref{in-02}),
   let $\tau$ be a positive real number and $u\not\equiv 0$ be any function in $W^{1,2}(\Sigma,g)$ satisfying
  \be\label{leq1}\|u\|_{1,\tau}^2=\int_\Sigma(|\nabla_g u|^2+\tau u^2)dv_g\leq 1.\ee
  By the Young inequality, one has for any $\epsilon>0$,
  \be\label{Young-ineq}u^2\leq (1+\epsilon)(u-\overline{u})^2+(1+\f{1}{4\epsilon})\overline{u}^2,\ee
   where
  \be\label{mean}\overline{u}=\f{1}{{\rm Area}(\Sigma)}\int_\Sigma udv_g.\ee
  Since $\int_\Sigma|\nabla_gu|^2dv_g\leq1-\tau\int_\Sigma u^2dv_g$ by (\ref{leq1}), we have
  $$\|u-\overline{u}\|_{1,\alpha}^2=\int_\Sigma|\nabla_gu|^2dv_g-\alpha\int_\Sigma(u-\overline{u})^2dv_g\leq 1-\tau\int_\Sigma u^2dv_g
  -\alpha\int_\Sigma(u-\overline{u})^2dv_g.$$
  Thus $\|u-\overline{u}\|_{1,\alpha}^2<1$ for $0\leq\alpha<\lambda_{\rm N}(\Sigma)$.  As a consequence,
   we can choose $\epsilon>0$ verifying
  $1+\epsilon={1}/{\|u-\overline{u}\|_{1,\alpha}^2}$.
  This leads to
  \be\label{55}\f{1}{\epsilon}=\f{\|u-\overline{u}\|_{1,\alpha}^2}{1-\|u-\overline{u}\|_{1,\alpha}^2}\leq \f{1-\tau\int_\Sigma u^2dv_g
  -\alpha\int_\Sigma(u-\overline{u})^2dv_g}
  {\tau\int_\Sigma u^2dv_g}\leq \f{1}{\tau\int_\Sigma u^2dv_g}.\ee
  Combining (\ref{mean}) and (\ref{55}), we get
  \be\label{u-bar}|\overline{u}|\leq\f{1}{{\rm Area}(\Sigma)}\int_\Sigma u^2dv_g\ee
  and thus
  \be\label{ep-1}\f{1}{4\epsilon}\overline{u}^2\leq \f{1}{4\tau{\rm Area}(\Sigma)}.\ee
  In view of (\ref{Young-ineq}), (\ref{u-bar}) and (\ref{ep-1}), we obtain
  \bna
  \exp(2\pi u^2)\leq \exp\le(2\pi\f{(u-\overline{u})^2}{\|u-\overline{u}\|_{1,\alpha}^2}\ri)
  \exp\le(2\pi\le(1+\f{1}{4\epsilon}\ri)\overline{u}^2\ri)
  \leq C\exp\le(2\pi\f{(u-\overline{u})^2}{\|u-\overline{u}\|_{1,\alpha}^2}\ri)
  \ena
  for some uniform constant $C$. Hence by applying (\ref{ineq-1}) of Theorem \ref{thm1},
  we conclude
  $$\int_\Sigma \exp(2\pi u^2)dv_g\leq C$$
  for some uniform constant $C$. Therefore (\ref{in-02}) follows immediately. \\

   {\it The inequality (\ref{in-02}) implies the inequality (\ref{in-01}).}

  Assume that (\ref{in-02}) holds. To prove (\ref{in-01}), we use the method of blow-up analysis. Suppose that
  (\ref{in-01}) does not hold. By (\ref{ineq-2}) for any $\gamma<2\pi$, we
  let $u_k$ be as in Lemma \ref{subcritical}. Then we must have
  \be\label{contrary}
  \lim_{k\ra\infty}\int_\Sigma \exp(\gamma_ku_k^2)dv_g=+\infty.
  \ee
  As before we assume with no loss of generality,  $c_k=u_k(x_k)=\max_{\Sigma}|u_k|$ and $x_k\ra x_0$ as $k\ra\infty$.
  Then the assumption (\ref{contrary}) implies
  that $c_k\ra+\infty$ as $k\ra\infty$.
  By Lemma \ref{Green},
  $c_ku_k$ converges to $G_{\alpha,x_0}$ strongly in $L^2(\Sigma,g)$.
  Since $\|u_k\|_{1,\alpha}=1$,
  \bna\gamma_ku_k^2&=&\gamma_k\f{u_k^2}{\|u_k\|_{1,\tau}^2}\|u_k\|_{1,\tau}^2\\
  &=&\gamma_k\f{u_k^2}{\|u_k\|_{1,\tau}^2}\le(\|u_k\|_{1,\alpha}^2+\alpha\int_\Sigma u_k^2dv_g+\tau\int_\Sigma u_k^2dv_g\ri)\\
  &=&\gamma_k\f{u_k^2}{\|u_k\|_{1,\tau}^2}+\f{\alpha\gamma_k}{\|u_k\|_{1,\tau}^2}u_k^2\int_\Sigma u_k^2dv_g+\f{\tau\gamma_k}
  {\|u_k\|_{1,\tau}^2}u_k^2\int_\Sigma u_k^2dv_g,
  \ena
  $\gamma_k=2\pi-1/k$, $\|u_k\|_{1,\tau}=1+o_k(1)$, $\alpha<\lambda_{\rm N}(\Sigma)$ and
  $$u_k^2\int_\Sigma u_k^2dv_g\leq \int_\Sigma c_k^2u_k^2dv_g=\int_\Sigma G_{\alpha,x_0}^2dv_g+o_k(1),$$
  we conclude
  $$\exp(\gamma_k u_k^2)\leq C\exp(\gamma_k u_k^2/\|u_k\|_{1,\tau}^2)$$
  for some uniform constant  $C$. It follows from (\ref{in-02}) that
  \bna\lim_{k\ra\infty}\int_\Sigma \exp(\gamma_k u_k^2)dv_g\leq C\lim_{k\ra\infty}\int_\Sigma \exp(\gamma_k u_k^2/\|u_k\|_{1,\tau}^2)dv_g\leq C\ena
  for some constant $C$. This contradicts (\ref{contrary}) and leads to (\ref{in-01}) immediately. $\hfill\Box$

\bigskip

{\bf Acknowledgement}. This work is partly supported by the National Science Foundation of China (Grant No.
  11761131002).


\begin{thebibliography}{00}

\bibitem{Adams} D. Adams, A sharp inequality of J. Moser for higher order derivatives, Ann. Math. 128 (1988) 385-398.

\bibitem{Adi-Stru} A. Adimurthi, M. Struwe, Global compactness
properties of semilinear elliptic equation with critical exponential
growth, J. Funct. Anal. 175 (2000) 125-167.

\bibitem{Adi-Yang} Adimurthi, Y. Yang, An interpolation of Hardy inequality
and Trudinger-Moser inequality in $\mathbb{R}^N$ and its
applications, Int. Math. Res. Notices 13 (2010)
2394-2426.

\bibitem{Aubin} T. Aubin, Nonlinear analysis on manifolds, Springer, 1982.

\bibitem{Bers} L. Bers, Riemann surfaces, Courant Institure Lecture Notes, 1957-58.

\bibitem{Bonheure} D. Bonheure, E. Serra, M. Tarallo, Symmetry of extremal functions in Moser-Trudinger inequalities and a H\'enon type problem
in dimension two, Adv. Differential Equations 13 (2008) 105-138.

\bibitem{Chang-Yang} A. Chang, P. Yang, Conformal deformation of metrics on $S^2$, J. Differential Geometry 27 (1988) 259-296.

\bibitem{Deng} S. Deng, New solutions for critical Neumann problems in $\mathbb{R}^2$, Adv. Nonlinear Anal. 8 (2019) 615-644.

\bibitem{Deng-Musso} S. Deng, M. Musso, Critical points of the Trudinger-Moser trace functional with high energy levels,
Ann. Inst. H. Poincare Anal. Non Lineaire 32 (2015) 59-95.

\bibitem{DJLW1} W. Ding, J. Jost, J. Li, G. Wang, An analysis of the two vetex case in the Chern-Simons Higgs model, Calc. Var. Partial
Differential Equation 7 (1998) 87-97.

\bibitem{DJLW2} W. Ding, J. Jost, J. Li, G. Wang, Self duality equations for Ginzburg-Landau and Seiberg-Witten type functionals with 6th
order potentials, Comm. Math. Phys. 217 (2001) 383-407.


\bibitem{DJLW} W. Ding, J. Jost, J. Li, G. Wang, The differential equation $\Delta u=8\pi-8\pi h e^u$
 on a compact Riemann surface, Asian
J. Math. 1 (1997) 230-248.

\bibitem{do-Yang} J. M. do \'O, Y. Yang, A quasi-linear elliptic equation with critical growth on compact Riemannian manifold without
boundary, Ann. Global Anal. Geom. 38 (2010) 317-334.

\bibitem{DRW} O. Druet, F. Robert, J. Wei, The Lin-Ni's problem for mean convex domains, Mem. Amer. Math. Soc. 218 (2012) 1027.

\bibitem{Fontana} L. Fontana, Sharp bordline Sobolev inequalities on compact Riemannian manifolds, Comment. Math. Helv. 68
(1993) 415-454.

\bibitem{GT} D. Gilbarg, N. Trudinger, Elliptic partial differential equations of second order, Springer, 2001.

\bibitem{Lan-Li} X. Lan, J. Li, Asymptotic behavior of the Chern-Simons Higgs 6th theory, Comm. Partial Differential Equations 32
(2007) 1473-1492.

\bibitem{LiJPDE} Y. Li, Moser-Trudinger inequality on compact
Riemannian manifolds of dimension two, J. Partial Differ.
Equ. 14 (2001) 163-192.

\bibitem{Li-Liu}  Y. Li, P. Liu, Moser-Trudinger inequality on the
boundary of compact Riemannian surface, Math. Z. 250 (2005)
363-386.

\bibitem{Liu} P. Liu, A Moser-Trudinger type inequality and blow up analysis
on compact Riemannian surface, Doctoral thesis, Max-Plank Institute, Germany, 2005.

\bibitem{Lu-Yang} G. Lu, Y. Yang, A sharpened Moser-Pohozaev-Trudinger inequality with mean value zero in $\mathbb{R}^2$, Nonlinear
Anal. 70 (2009) 2992-3001.

\bibitem{Moser} J. Moser, A sharp form of an inequality by N. Trudinger, Indiana Univ. Math. J. 20 (1971) 1077-1092.

\bibitem{Nguyen} Q. Ng\^o, V. Nguyen, An improved Moser-Trudinger inequality involving the first non-zero Neumann eigenvalue with mean value zero
in $\mathbb{R}^2$, arXiv: 1702.08883.

\bibitem{Nguyen-2} V. Nguyen, A sharp Adams inequality in dimension four and its extremal functions, arXiv: 1701.08249.


\bibitem{Peetre} J. Peetre, Espaces d'interpolation et theoreme de Soboleff, Ann. Inst. Fourier (Grenoble) 16 (1966) 279-317.

\bibitem{Pohozaev} S. Pohozaev, The Sobolev embedding in the special case
$pl=n$, Proceedings of the technical scientific conference on
advances of scientific reseach 1964-1965, Mathematics sections,
158-170, Moscov. Energet. Inst., Moscow, 1965.

\bibitem{Pom} C. Pommerenke, Boundary behavior of conformal maps, Springer, 1992.

\bibitem{Tintarev} C. Tintarev, Trudinger-Moser inequality with remainder terms, J. Funct. Anal. 266 (2014) 55-66.

\bibitem{Trudinger} N. Trudinger, On embeddings into Orlicz spaces and
some applications, J. Math. Mech. 17 (1967) 473-484.


\bibitem{WangM1} M. Wang, The self-dual Chern-Simons Higgs equation on a compact Riemann surface with boundary, Internat. J. Math. 21
(2010) 67-76.

\bibitem{WangM-2} M. Wang, The asymptotic behavior of Chern-Simons Higgs model on a compact Riemann surface with boundary, Acta Math. Sin.
(Engl. Ser.) 28 (2012) 145-170.

\bibitem{2006} Y. Yang, Extremal functions for Moser-Trudinger inequalities on 2-dimensional compact Riemannian manifolds with boundary,
Internat. J. Math. 17 (2006) 313-330.

\bibitem{Yang-PJM} Y. Yang, Moser-Trudinger trace inequalities on a compact Riemannian surface with boundary, Pacific J. Math.
227 (2006) 177-200.


\bibitem{Yang-MZ} Y. Yang, A sharp form of trace Moser-Trudinger inequality on compact Riemannian surface with boundary,
Math. Z. 255 (2007) 373-392.

\bibitem{Yang-Tran} Y. Yang, A sharp form of the Moser-Trudinger inequality on a compact Riemannian surface,
Trans. Amer. Math. Soc. 359 (2007) 5761-5776.

\bibitem{Yang-JDE-15} Y. Yang, Extremal functions for Trudinger-Moser inequalities of Adimurthi-Druet type in dimension two,
J. Differential Equations 258 (2015) 3161-3193.

\bibitem{Yudovich} V. Yudovich, Some estimates connected with integral operators and with solutions of elliptic equations,
Sov. Math. Docl. 2 (1961) 746-749.

\bibitem{ZhouChunqin} T. Zhang, C. Zhou, Asymptotical behaviors for Neumann boundary problem with singular data, Acta Math. Sin.
(Engl. Ser.) 35 (2019) 463-480.

\bibitem{ZhuX} X. Zhu, Solutions for Toda system on Riemann surface with boundary, Acta Math. Sin. (Engl. Ser.) 27 (2011) 1501-1520.
\end{thebibliography}
\end{document}